\documentclass[11pt]{amsart}

\usepackage{amsmath,amsthm,amssymb,url}
\usepackage{verbatim}

\usepackage[all]{xy}

\DeclareMathOperator{\Spec}{Spec}
\DeclareMathOperator{\Hom}{Hom}
\DeclareMathOperator{\Id}{Id}

\DeclareMathOperator{\Def}{Def}

\DeclareMathOperator{\Der}{Der}
\DeclareMathOperator{\Coder}{Coder}

\DeclareMathOperator{\MC}{MC}
\DeclareMathOperator{\KS}{KS}

\DeclareMathOperator{\tr}{Tr}

\newcommand{\C}{\mathbb{C}}
\newcommand{\K}{\mathbb{K}}
\newcommand{\Z}{\mathbb{Z}}
\newcommand{\N}{\mathbb{N}}

\newcommand{\D}{\mathcal{D}}
\newcommand{\T}{\mathcal{T}}
\newcommand{\A}{\mathcal{A}}
\newcommand{\deltabar}{\bar{\partial}}
\newcommand{\Eps}{\mathcal{E}}
\newcommand{\U}{\mathcal{U}}
\newcommand{\Oh}{\mathcal{O}}
\renewcommand{\l}{\mathfrak{l}}

\newcommand{\contr}{{\mspace{1mu}\lrcorner\mspace{1.5mu}}}

\newtheorem{theorem}{Theorem}[section]
\newtheorem{lemma}[theorem]{Lemma}
\newtheorem{proposition}[theorem]{Proposition}
\newtheorem{corollary}[theorem]{Corollary}

\theoremstyle{definition}
\newtheorem{definition}[theorem]{Definition}
\newtheorem{example}[theorem]{Example}
\newtheorem{remark}[theorem]{Remark}

\newenvironment{acknowledgement}{\par\addvspace{17pt}\small\rm
\trivlist\item[\hskip\labelsep{\it Acknowledgement.}]}
{\endtrivlist\addvspace{6pt}}

\begin{document}
\title{Infinitesimal deformations of Hitchin pairs and Hitchin map}
\author{Elena Martinengo}
\email{martinengo@mat.uniroma1.it}
\urladdr{www.mat.uniroma1.it/people/}

\begin{abstract}
We identify dglas that control infinitesimal deformations of the pairs (manifold, Higgs bundle) and of Hitchin pairs. As a consequence, we recover known descriptions of first order deformations and we refine known results on obstructions.   
Secondly we prove that the Hitchin map is induced by a natural $L_\infty$-morphism and, by standard facts about $L_\infty$-algebras, we obtain new conditions on obstructions to deform Hitchin pairs.
\end{abstract}

\subjclass{17B70, 14D15, 14D20}
\keywords{Hitchin pairs, Differential graded Lie algebras, $L_\infty$-algebras}

\maketitle

\section*{Introduction}

The interest on Higgs bundles started twenty years ago with Nigel Hitchin and Carlos Simpson's studies.
The rich structure of Higgs bundles manage they play a role in many different mathematical areas.
A Higgs bundle on a complex manifold $X$ is a pair $(E,\theta)$, where $E$ is a holomorphic vector bundle on $X$ and $\theta \in H^0(X, \Eps nd(E)\otimes \Omega_X^1)$, such that $\theta \wedge \theta=0$. A generalization of Higgs bundles are Hitchin pairs $(E, L, \theta)$, where $L$ is a bundle on $X$ and now $\theta \in H^0(X, \Eps nd(E)\otimes L)$.

\smallskip

In the present work we study infinitesimal deformations of Higgs bundles and Hitchin pairs and give a description of the Hitchin map as a morphism of deformation theories; we use differential graded Lie algebras to analyse these  deformations. 

The philosophy underlying this approach, originating in the works of Quillen,
Deligne, Drinfeld and Kontsevich, is that, in characteristic zero, every deformation problem
is governed by a differential graded Lie algebra, via the deformation functor associated
to it, given by solutions of Maurer-Cartan equation modulo gauge action.
Dglas tecniques allow to preserve a lot
of informations on the deformation problem, which are lost with classical methods, and some classical results can be obtained as easy consequences of
definitions and formal constructions. 

One of our main goals is to find out dglas that govern infinitesimal deformations of a Higgs bundle, of a pair (manifold, Higgs bundle) and of a Hitchin pair. Biswas and Ramanan in \cite{Biswas, Biswas-Ramanan} introduced complexes of sheaves to study these deformations: for a Higgs bundle $(E, \theta)$ the complex to be considered is
\[ 0 \to \Eps nd(E) \stackrel{[-,\theta]}{\longrightarrow} \Eps nd (E) \otimes \Omega_X^1 \stackrel{[-,\theta]}{\longrightarrow} \Eps nd (E) \otimes \Omega_X^2\stackrel{[-,\theta]}{\longrightarrow}  \Eps nd (E) \otimes \Omega_X^3\to \ldots,\]
where the differential is defined using the composition of endomorphisms and the wedge product of forms, and for a pair $(X, (E, \theta))$ the complex is
\[ 0 \to \D^1(E) \stackrel{[-,\theta]}{\longrightarrow} \Eps nd (E) \otimes \Omega_X^1 \stackrel{[-,\theta]}{\longrightarrow} \Eps nd (E) \otimes \Omega_X^2\stackrel{[-,\theta]}{\longrightarrow}  \Eps nd (E) \otimes \Omega_X^3\to \ldots ,\] 
where the first differential is given by the sum of the above multiplication by $\theta$ and the action of differential operators on forms via the Lie derivative. 
We obtain the dglas that govern infinitesimal deformations of a Higgs bundles and of a pair (manifold, Higgs bundle) defining  dgla structures on the total complexes of the Dolbeault resolutions of the above complexes of sheaves. From the identification of these dglas, we are immediately able to recover and refine Biswas and Ramanan results about tangent and obstruction spaces. 
We prove that the spaces of first order deformations of a Higgs bundle $(E, \theta)$ and of a pair $(X,(E,\theta))$ are naturally isomorphic to the first hypercohomology spaces respectively of the first and of the second complex of sheaves, 
while obstructions are contained in the second hypercohomology spaces of them (Theorem \ref{Teo.dgla controlla higgs bundle} and Corollary \ref{Cor.dgla Higgs bundle}). For the general case of a Hitchin pair $(E,L, \theta)$, a similar result holds, substituting, in the first of the above complexes, the sheaf $\Omega_X^*$ with  $\bigwedge^*L$ (Theorem \ref{Teo.dgla L Higgs bundle}).

\smallskip

Secondly we concentrate our attention on the study of the Hitchin map from deformations point of view. It is defined as:
\[ H: {\mathcal M} \to \bigoplus_{k=1}^{r} H^0(X, \bigodot^k L), \qquad  H(E, \theta) = (\tr(\theta), \ldots, \tr(\theta^r)),\]
from the moduli space of Hitchin pairs on $X$ to the space of global sections of the vector bundles $\bigodot^k L$. 
This map was introduced for Higgs bundles by Simpson in \cite{Simpson 2}, as a generalization of the determinant map that Hitchin studied on curves in \cite{Hitchin}; a version of this map for Hitchin pairs can be found in \cite{Biswas-Gothen-Logares}. 

To study the Hitchin map in terms of deformation theory, the dglas approach is not convenient and we use the more powerfull tool of $L_\infty$-algebras.
Theory of deformations via differential graded Lie algebras and via $L_\infty$-algebras is based on the principle that the local study of a moduli space is encoded by a dgla or an $L_\infty$-algebra opportunely choosen. Then every natural morphism between moduli spaces is induced by a morphism between the associated dglas or $L_\infty$-algebras.

If we intend the Hitchin map as a morphism of moduli spaces, it is natural to expect it is induced by a morphism between the algebraic objects associated to them. Since the Hitchin map is not even linear, we can not expect to obtain it from a dglas morphism. 
In Propositions \ref{Prop.h L-infty morphism} and \ref{Prop.h Hitchin map}, we explicitate an $L_\infty$-morphism $h$ that induces the Hitchin map. 
As a direct consequence of this $L_\infty$-description and of $L_\infty$-tecniques, in Corollary  \ref{Cor.ostruzione e mappa Hitchin}, we obtain that obstructions to deform Hitchin pairs are contained in the kernel of the map induced at second cohomology level by the linear part of $h$.  
It is defined on the second hypercohomology space of the complex of sheaves
\[ 0 \to \Eps nd(E) \stackrel{[-,\theta]}{\longrightarrow} \Eps nd (E) \otimes L \stackrel{[-,\theta]}{\longrightarrow} \Eps nd (E) \otimes (L\wedge L)\stackrel{[-,\theta]}{\longrightarrow}  \Eps nd (E) \otimes (L\wedge L \wedge L)\to \ldots,\] 
where the differential is the multiplication by $\theta \in H^0(X, \Eps nd(E) \otimes L)$, defined using the composition of endomorphisms and the wedge product on forms and on $L$. 
Taking the Dolbeault resolution, this map associates to every element $\omega\otimes f \in A_X^{0,1}(\Eps nd (E) \otimes L)$ the class in $\bigoplus_{k=1}^{r} H^1(X, \bigodot^k L)$ of multiples of $\omega \otimes \tr(f\theta^{k-1})$. 

\medskip

The paper is organized as follows. In the first section, we introduce some basic definitions of dglas theory, that we use in all the paper and in the second section to analyse deeply deformations of a complex manifold and of its holomorphic forms. The third section is devoted to deformations of Higgs bundles and of pairs (manifold, Higgs bundle), with the statements of our main results: the identification of dglas that govern these deformations and the description of the tangent and obstruction spaces. These results are proved in all details in section \ref{Sect.proof of main theorem}.   
In the next section, we study deformations of Hitchin pairs, using the same tecniques and obtaining similar results.
In the sixth section, we introduce some basic definitions and tools of $L_\infty$-algebras theory, that are essential to the study of the Hitchin map, we do in the last section. There we give a deformation theoretic interpretation of the Hitchin map obtaining it from an $L_\infty$-morphism and we deduce a condition on obstructions of Hitchin pairs.

\begin{acknowledgement} I thank Marco Manetti for having introduced me to the problem and for all his precious advices  and suggestions on the preparation of this paper. I also thank Domenico Fiorenza for some stimulating discussions on subjects related to the paper.  
\end{acknowledgement}

\section{Deformation theory via dglas} \label{sec.dgla}
In this section we introduce some tools of deformation theory via differential graded Lie algebras: we give the basic definitions and we analyse classical examples. The main references we follow are \cite{Man Pisa, Man Roma}.  
\begin{definition}
A \emph{differential graded Lie algebra}, dgla, is the data $(L,d,[\ ,\ ])$, where $L=\bigoplus_{i\in \mathbb {Z}} L^i$ is a $\Z$-graded vector space over $\C$, $d:L^i \rightarrow L^{i+1}$ is a linear map, such that $d \circ d=0$, and $[\ ,\ ]:L^i \times L^j \rightarrow L^{i+j}$ is a bilinear map, such that:
\begin{enumerate}
\item[-] $[\ ,\ ]$ is graded skewsymmetric, i.e. $[a,b]=-(-1)^{\deg a\deg b}[b,a]$,  
\item[-] $[\ ,\ \ ]$ verifies the graded Jacoby identity, i.e. $[a,[b,c]]=[[a,b],c]+(-1)^{\deg a\deg b}[b,[a,c]]$,
\item[-] $[\ ,\ ]$ and $d$ verify the graded Leibniz's rule, i.e. $d[a,b]=[da,b]+(-1)^{\deg a}[a,db]$, 
\end{enumerate}
for every $a, b$ and $c$ homogeneous.
\end{definition}
\begin{definition}
Let $(L,d_L, [\ ,\ ]_L)$ and $(M, d_M, [\ ,\ ]_M)$ be two dglas, a \emph{morphism of dglas} $\phi:L\to M$ is a degree zero linear morphism that commutes with the brackets and the differentials.  
\end{definition}

\noindent We now introduce some usefull examples of differential graded Lie algebras. 
\begin{example} \label{Ex.dgla prod tensoriale}
Let $(A=\oplus_i A_i, d_A)$ be a differential graded commutative $\C$-algebra and let $(L=\oplus_i L^i, d_L, [\ , \ ])$ be a dgla, then the tensor product $A\otimes_\C L$ has a natural structure of dgla by setting:
\[ (A\otimes_\C L)^n= \bigoplus_{p+q=n} A_p \otimes_\C L^q, \]
\[  d(a\otimes x)= d_A a \otimes x + (-1)^{\deg a} a \otimes d_L x, \qquad [a\otimes x, b\otimes y]= (-1)^{\deg b \cdot \deg x} a\cdot b \otimes [x,y], \]
for all $a, b \in A$ and $x,y \in L$. 
\end{example}

\begin{example} \label{Ex.Hom dgla}
Let $(V=\bigoplus_{i\in \Z} V^i, d)$ be a differential $\Z$-graded $\C$-vector space. Consider the $\Z$-graded $\C$-vector space $\Hom(V,V) = \bigoplus_{i\in \Z}  \Hom^i(V,V)$,  
where $\Hom^i(V,V) =\{ f: V \to V \mbox{\ \ linear} \mid f(V^n) \subset f(V^{n+i}) \mbox{\ \ for every\ } n\}$.
The bracket \[[f,g] = f \circ g  -(-1)^{\deg f \deg g}g\circ f\]
and the differential \[df = [d,f] = d\circ f - (-1)^{\deg f} f\circ d\]
make $\Hom(V,V)$ a differential graded Lie algebra.
Now let $L$ be a dgla, consider the subspace $\Der(L,L) \subset \Hom(L,L)$ of derivations of $L$, where a derivation $f:L\to L$ is a linear map, which satisfies the graded Leibniz's rule:
\[ f([a,b])= [f(a),b]+(-1)^{\deg f \cdot \deg a} [a, f(b)].\]
$\Der(L,L)$ with the above differential and bracket is a sub-dgla of $\Hom(L,L)$.   
\end{example}

\begin{example} \label{Ex.KS dgla}
Let $X$ be a compact complex manifold, let $\mathcal{T}_X$ be the holomorphic tangent bundle of $X$, let $\A_X^{p,q}(\T_X)=\A_X^{p,q}\otimes_{\Oh_X}\T_X$ be the sheaf of $(p,q)$-forms of $X$ with values on the tangent bundle. The \emph{Kodaira-Spencer} dgla of $X$ is 
\[ \KS_X = \bigoplus_{i\in \N} \Gamma(X, \A^{0,i}_X(\T_X))=\bigoplus_{i\in \N} A^{0,i}_X(\T_X), \]
the space of the global sections of the sheaf of $(0,i)$-forms of $X$ with values on the tangent bundle $\T_X$. 
The dgla structure on $\KS_X$ is given as follows. 
The differential on $\KS_X$ is the Dolbeault differential and the bracket is defined in local coordinates, $z_1, \ldots , z_n$, extending the standard bracket on $\A_X^{0,0}(\T_X)$ bilinearly with respect to the sheaf of the antiholomorphic differential forms:
\[
[f d\bar{z}_I \frac{\partial}{\partial z_i}, g d\bar{z}_J \frac{\partial}{\partial z_j}]  =  (f \frac{\partial g}{\partial z_i} \frac{\partial}{\partial z_j} -  g\frac{\partial f}{\partial z_j} \frac{\partial}{\partial z_i}) d\bar{z}_I\wedge d\bar{z}_J,   \qquad \forall f, g \in \A^{0,0}_X. \]
\end{example}

\begin{example} \label{Ex.End(E),D(E) dgla}
Let $X$ be a complex manifold, let $E$ be a locally free sheaf of $\Oh_X$-modules on $X$ and $\Eps nd (E)$ the sheaf of endomorphisms of $E$. 
We indicate with $\D iff^1 (E)$ the sheaf of differential operators of degree $\leq 1$ on sections of $E$ and with $\D^1(E)\subset \D iff^1 (E)$ the subsheaf of operators with scalar principal symbol $\sigma$, defined by the following exact sequence of sheaves:
\[0 \to \Eps nd (E) \to  \D^1(E) \stackrel{\sigma}{\rightarrow} \mathcal{T}_X \to 0. \]  
Consider the graded vector spaces: 
\[ \bigoplus_{i\in \N}  A^{0,i}_X(\Eps nd(E))\quad \mbox{and}\quad \bigoplus_{i\in \N}  A^{0,i}_X(\D^1(E)), \]
dgla-structures are defined on them taking the following differential and bracket:
\[ d(\omega\otimes P)= \deltabar \omega \otimes P, \]
\[ [\omega\otimes P, \eta\otimes Q]= \omega\wedge \eta \otimes[P,Q] + \omega\wedge \l_{\sigma(P)}(\eta) \otimes Q -(-1)^{\deg\omega \cdot \deg \eta}   \eta \wedge \l_{\sigma(Q)}(\omega) \otimes P, \]
for all $\omega, \eta \in A_X^{0,*}$ and $P,Q\in \D^1(E)$, the symbol $\l$ indicates the Lie derivative (Section \ref{Sect. def holomorphic forms}). Note that, if $P,Q\in \Eps nd(E)$, their symbol is zero and we obtain the classical dgla structure on $A_X^{0,*}(\Eps nd (E))$ and that, if $P,Q \in A_X^{0,*}(\T_X)$, we recover the dgla structure of Example \ref{Ex.KS dgla}, thus we have the exact sequence of dglas:
\[ 0 \to A_X^{0,*}(\Eps nd (E)) \to  A_X^{0,*}(\D^1(E)) \stackrel{\sigma}{\rightarrow} A_X^{0,*}(\mathcal{T}_X) \to 0.\]
\end{example}

\medskip
\noindent The differential graded Lie algebras approach to deformation theory is based on the following definition of the deformation functor associated to a dgla. We indicate with $\bf{Art}_\C$ the category whose objects are local Artinian $\C$-algebras and whose arrows are local morphisms of $\C$-algebras and with $\bf{Set}$ the category of sets.

\begin{definition}
Let $L$ be a dgla, the \emph{deformation functor} associated to $L$ is the functor $\Def_L: \bf{Art}_{\C}\to \bf{Set}$, given for all $A\in \bf{Art}_{\C}$, by:
$$\Def_L(A)=\frac{\MC_L(A)}{\sim_{\textrm{gauge}}},$$ 
where:  $\qquad \qquad \qquad \displaystyle \MC_L(A)=\left\{x\in L^1\otimes \mathfrak{m}_A \mid dx+\frac{1}{2}[x,x]=0\right\}$ \\
and the gauge action is the action of $\exp (L^0\otimes \mathfrak{m}_A)$ on $\MC_L(A)$, given by:
$$e^a * x= x+\sum_{n=0}^{+\infty} \frac{([a,-])^n}{(n+1)!}([a,x]-da).$$
\end{definition} 
It is evident that a morphism of dglas $f:L\to M$ induces a natural transformation of functors $\MC(f): \MC_L\to \MC_M$. It is compatible with the gauge action and induces a natural transformation of deformation functors $\Def(f): \Def_L \to \Def_M$.

\begin{remark} \label{rm.def dgla}
Let $L$ be a dgla and $\Def_L$ the deformation functor associated to $L$.  
It can be proved that the tangent space to $\Def_L$ is the first cohomology space of $L$, $H^1(L)$, and that obstructions are naturally contained in $H^2(L)$.\\
It is a classical calculation to observe that,
if $f: L\to M$ is a morphism of dglas, then the linear maps $H^1(f) : H^1(L)\to H^1(M)$ and $H^2(f) : H^2(L)\to H^2(M)$ are morphisms of tangent spaces and of obstruction spaces, respectively, compatible with the morphism $\Def(f): \Def_L \to \Def_M$.
\end{remark}

\begin{remark} \label{rm.mc vs gauge banale}
It is easy to verify that a deformation functor is trivial if and only if its tangent space is trivial. 
This allow to prove that, for dglas of type $A_X^{0,*}({\mathcal F})$, with $\mathcal F$ quasi-coerent sheaf on a complex manifold $X$, every Maurer-Cartan solution is locally gauge equivalent to zero. Indeed, taking an affine open cover $\U=\{U_i\}_i$ of $X$, $H^1(U_i, {\mathcal F})=0$ is the tangent space to the functor $\Def_{A^{0,*}_{U_i}({\mathcal F})}$ which is trivial.   
\end{remark}

\noindent Let $\mathcal X$ be a geometric object, for example a manifold or a sheaf, and let $\Def_{\mathcal X}$ be the functor of infinitesimal deformations of $\mathcal X$, i.e. the functor:
\[  \Def_{\mathcal X}: \bf{Art}_\C \to \bf{Set}, \]
that associates to every local Artinian $\C$-algebra the set of isomorphism classes of deformations of $\mathcal X$ over it.
If there exists a dgla $L$, such that $\Def_{\mathcal X}$ is isomorphic to the deformation functor associated to $L$, we say that the dgla $L$ governs deformations of $\mathcal X$. 
Let's state two well known examples of this situation:

\begin{example} \label{ex.def varieta}
Let $X$ be a complex manifold and let $\KS_X$ be the Kodaira-Spencer dgla of it. We indicate with $\Def_X$ the functor of infinitesimal deformations of $X$ and with $\Def_{\KS_X}$ the deformation functor associated to the dgla $\KS_X$. The natural transformation of functors:
$$\begin{array}{rrllr}
\Phi: & \Def_{\KS_X}(A)& \longrightarrow &\Def_{X}(A), & \qquad \forall A \in \bf{A}rt_{\C},  \\
& x& \longrightarrow & \ker(\deltabar+\mathfrak{l}_x)& 
\end{array}$$
where $\l_x$ is the holomorphic Lie derivative (Section \ref{Sect. def holomorphic forms}), is an isomorphism. In Section \ref{Sect. def holomorphic forms}, we will analyse in more details deformations of a complex manifold and the link with the Kodaira-Spencer dgla. In proofs of Propositions \ref{prop.def mappa} and \ref{prop.mappa sui def}, we will relate the principal steps in the construction of $\Phi$.
For a complete study of it, its construction and the proof that it is an isomorphism, see \cite{Dona.Tesi}.
\end{example}

\begin{example} \label{ex.def (X,E)}
Let $X$ be a complex manifold, let $E$ be a locally free sheaf of $\Oh_X$-modules on $X$ and let $A_X^{0,*}(\D^1(E))$ be the dgla of the $(0,*)$-forms with values in the first order differential operators with scalar principal symbol. We indicate with $\Def_{(X,E)}$ the functor of infinitesimal deformations of the pair $(X,E)$ and with $\Def_{A_X^{0,*}(\D^1(E))}$ the deformation functor associated to the dgla $A_X^{0,*}(\D^1(E))$.
The natural transformation of functors:
$$\begin{array}{rrllr}
\Psi: &\Def_{A_X^{0,*}(D^1(\Eps))}(A)& \longrightarrow &\Def_{(X,\Eps)}(A), & \qquad \forall A \in \bf{A}rt_{\C},  \\
& x& \longrightarrow & (\ker(\deltabar+\l_{\sigma(x)}),\ker(\deltabar+x)) &
\end{array}$$
is an isomorphism. In proofs of Propositions \ref{prop.def mappa} and \ref{prop.mappa sui def}, 
we will relate the principal steps in the construction of it. 
For a complete study of $\Psi$ see \cite{Tesi.mia}, where it is defined and it is proved to be an isomorphism.
\end{example}

\section{Infinitesimal deformations of holomorphic forms} \label{Sect. def holomorphic forms}
In this section we analyse deformations of a complex manifold and the consequent deformations of the sheaves of holomorphic functions and of holomorphic forms. We explain the link between deformations of a complex manifold and the Kodaira-Spencer dgla of it. We mainly follow \cite{Fio-Man periodi, Fio-Man periodi generalizzati}.

\medskip

Start with a differentiable manifold $X$, 
recall that an \emph{almost complex structure} on $X$ can be seen as a subsheaf $
{\mathcal V}\subset \A^1_X$ of locally free $\A^0_X$-modules, such that ${\mathcal V}\oplus \overline{{\mathcal V}} = \A^1_X$. 
Obviously, if $X$ is a complex manifold, the decomposition $\A^1_X=\A^{0,1}_X\oplus \overline{\A^{0,1}_X}$ define an almost complex structure on $X$. 
An almost complex structure on a differentiable manifold $X$ is called \emph{integrable} if there exist a structure of complex manifold on $X$ that induces it. 
An integrable almost complex structure is called a \emph{complex structure}.

Frobenius and Newlander-Nirenberg's Theorems (\cite[Ch. 2]{Voisin}) concern conditions under which an almost complex structure is integrable. 

A deformation of a complex manifold $X$ can be seen as a deformation of its complex structure, infact, by Ehresmann's Theorem (\cite[Theorem 2.4]{Kodaira}, \cite[Theorem 9.3]{Voisin}), deformations of a complex manifold are diffeomorphic.
Now we study deformations of a complex manifold from this point of view. 

We start with some definitions.
In general, for any   vector space $V$  and linear functional
$\alpha:V \to \C$, the \emph{contraction operator} is defined as 
$$
\alpha \contr: \bigwedge^k V \to \bigwedge^{k-1}V,
$$
$$
\alpha\contr (v_1 \wedge \ldots \wedge v_k)=\sum_{i=1}^k
(-1)^{i-1} \alpha(v_i)(v_1 \wedge \ldots \wedge \hat{v}_i \wedge
\ldots \wedge v_k);
$$
it is a derivation of degree $-1$ of the graded algebra
$(\bigwedge^k V ,\wedge)$.
The contraction  $\contr$ of differential
forms with vector fields defines an injective morphisms of
sheaves
\[ \begin{array}{llll}
{\boldsymbol{i}}:&  \A^{0,*}_X(\T_X) & \longrightarrow & \mathcal{D}er^*(\A_X^{*},\A_X^{*})[-1]\\
                 &   \xi               &  \longmapsto    & {\boldsymbol{i}}_\xi    \end{array} \]   
where ${\boldsymbol{i}}_\xi(\omega)=\xi \contr \omega$, for all $\omega \in \A_X^{*}$.\\
There is an other action of the sheaf $\A_X^{0,*}(\T_X)$ on the sheaf $\A_X^{*}$ as a derivation, via the holomorphic \emph{Lie derivative}: 
\[ \begin{array}{llll}
\l:& \A^{0,*}_X(\T_X) & \longrightarrow & Der^*(\A_X^{*}, \A_X^{*})\\
   & \xi                  & \longmapsto     & \l_\xi=[\partial ,{\boldsymbol{i}}_\xi] \end{array}\]
given by: 
\[\l_\xi(\omega)=[\partial ,{\boldsymbol{i}}_\xi] (\omega) =\partial(\xi\contr \omega)+(-1)^{\deg\xi} \xi \contr \partial\omega, \qquad \mbox{for all} \ \omega \in
\A_X^{*}.\]
Observe that the holomorphic Lie derivative define an action of the sheaf $\T_X$ on the sheaf $\Omega^*_X$ of the holomorphic forms on $X$: 
\[ \mathfrak{l}: \T_X \to \mathcal{D}er^0(\Omega^{*}_X,\Omega^{*}_X),  \]
infact, for all $x\in \T_X$, the derivation $\mathfrak l_x$ of a holomorphic form gives as result a holomorphic form of the same degree.\\
For future use, write out some properties of the holomorphic Lie derivative. 
\begin{lemma} \label{lemma eq cartan homotopy}
For every $\xi, \eta \in \A_X^{0,*}(\T_X)$, the following equalities hold:
\[   \boldsymbol{i}_{d\xi} = -[\deltabar, \boldsymbol{i}_\xi], \qquad \boldsymbol{i}_{[\xi, \eta]} = [\boldsymbol{i}_\xi,[\partial, \boldsymbol{i}_\eta]] \qquad \mbox{and} \qquad [\boldsymbol{i}_\xi, \boldsymbol{i}_\eta]=0.\]
\end{lemma}
\begin{proof}
See \cite[Lemma 2.1]{Man Kahler}.
\end{proof}
\begin{lemma}
The holomorphic Lie derivative $\mathfrak{l}: A_X^{0,*}(\T_X)\to \Der^*(A^{0,*}_X,A^{0,*}_X)$ is a morphism of dglas.
\end{lemma}
\begin{proof} Compatibility with differentials is given using definitions and the Jacobi identity: 
\[ d(\l_\xi) = [\deltabar, \l_\xi] = [\deltabar, [\partial, {\boldsymbol{i}}_\xi]] = -[\partial, [ \deltabar, {\boldsymbol{i}}_\xi]] = [\partial,{\boldsymbol{i}}_{d\xi} ]= \l_{d\xi}, \]
for all $\xi \in  A_X^{0,*}(\T_X)$. Compatibility with brackets is similar.
\end{proof}
The above lemma implies that $\xi\in A_X^{0,1}(\T_X)$ satisfies the Maurer-Cartan equation if and only if $\deltabar + \l_\xi$ is a differential. Infact 
\[ (\deltabar +\l_\xi)^2 = \deltabar^2 + \deltabar \l_\xi + \l_\xi \deltabar + \l_\xi^2= [\deltabar, \l_\xi]+ \frac{1}{2}[\l_\xi,\l_\xi]= \l_{d\xi+ \frac{1}{2}[\xi,\xi]}.\] 

Now we are ready to study deformations of a complex manifold $X$. Fix $\xi \in A_X^{0,1}(\T_X)\otimes \mathfrak{m}_A$, with $A\in \bf{Art}_\C$.  
Define a deformation of almost complex structure associated to $\xi$:
\[ \A_\xi^{1,0} = \{ \omega \in \A_{X}^1 \otimes A \mid \pi_{0,1}(\omega) = \boldsymbol{i}_\xi (\pi_{1,0}(\omega))\},  \]
where $\pi_{0,1}: \A^1_X \to \A_X^{0,1}$ and $\pi_{1,0}: \A^1_X \to \A_X^{1,0}$ are the projections. 
Newlander-Nirenberg's Theorem assures that the almost complex structure $\A_{\xi}^{1,0}$ is integrable if and only if $\xi$ is a solution of the Maurer-Cartan equation.
Therefore every element $\xi \in \MC_{\KS_X}(A)$ define a deformation of complex structure of $X$.

Let $\xi \in A^{0,1}_X(\T_X) \otimes \mathfrak{m}_A$ be a Maurer-Cartan element and let $\A_{\xi}^{0,1}$ be the deformation of complex structure of $X$ associated to $\xi$ as before, this defines
\[ \A^{0,1}_\xi = \overline{\A_{\xi}^{1,0}}, \qquad \A^{p,q}_{\xi}=\bigwedge^p \A_{\xi}^{1,0} \otimes \bigwedge^q \A_{\xi}^{0,1}.\]
Moreover, by definition, the sheaf of $\xi$-holomorphic function is
\[ \Oh_\xi=\{ f\in \A_X^0 \otimes A \mid df\in \A_\xi^{1,0} \} = \{f\in \A_X^0 \mid \deltabar f =\xi \contr \partial f \} = \{f\in \A_X^0 \mid (\deltabar +\l_\xi) f=0\}.\]
To describe how the sheaf of holomorphic differential forms change with the deformation of complex structure, we need the following:
\begin{lemma} \label{Lemma.e iso}
Let $\xi \in A^{0,1}_X(\T_X)\otimes \mathfrak{m}_A$, the exponential 
\[ e^{\boldsymbol{i}_\xi}: \A^*_X \otimes A \to \A^*_X\otimes A  \]
is an isomorphism of graded algebras, such that $e^{\boldsymbol{i}_\xi}(\A_X^{1,0} \otimes A)= \A_\xi^{1,0}$.
\end{lemma}
\begin{proof} It is an easy application of definitions, we follow \cite[Lemma 11.2]{Fio-Man periodi}.The only fact to be proved is the equality $e^{\boldsymbol{i}_\xi}(\A_X^{1,0} \otimes A)= \A_\xi^{1,0}$. Observe that, if $\omega \in \A^1_X\otimes A$: \[e^{\boldsymbol{i}_\xi}(\omega) = \omega + \xi \contr \omega= \omega + \xi \contr \pi_{1,0}(\omega),\] 
then $\pi_{1,0}(e^{\boldsymbol{i}_\xi}(\omega)) = \pi_{1,0}(\omega)$ and $\pi_{0,1}(e^{\boldsymbol{i}_\xi}(\omega)) = \pi_{0,1}(\omega) + \xi \contr \pi_{1,0}(\omega).$
Therefore, using the above definitions, $e^{\boldsymbol{i}_\xi}(\omega)\in  \A_\xi^{1,0}$ if and only if $\pi_{0,1}(\omega)=0$.  
\end{proof}
\begin{lemma}Let $\xi \in A^{0,*}_X(\T_X) \otimes \mathfrak{m}_A$ be a solution of the Maurer-Cartan equation, then
\[ e^{-\boldsymbol{i}_\xi} d e^{\boldsymbol{i}_\xi} = d + e^{-\boldsymbol{i}_\xi}* 0= d +\l_\xi,  \]
in the dgla $\Der^*(A^*_X,A_X^*)\otimes A$. 
\end{lemma}
\begin{proof} A proof can be found in \cite[Corollary 12.2]{Fio-Man periodi}. Here we make explicit calculations.
The first equality is given by definition of gauge action in the dgla $\Der^*(A_X^*, A_X^*)$:
\[ e^{-\boldsymbol{i}_\xi}*0 = 0+\sum_{n=0}^{+\infty} \frac{([-\boldsymbol{i}_\xi,-])^n}{(n+1)!}\ ([-\boldsymbol{i}_\xi,0]+d\boldsymbol{i}_\xi) = 
\sum_{n=0}^{+\infty} \frac{([-\boldsymbol{i}_\xi,-])^n}{(n+1)!}\ [-\boldsymbol{i}_\xi,d]=\]\[= \sum_{n=1}^{+\infty} \frac{([-\boldsymbol{i}_\xi,-])^n}{n!}\ d=
\sum_{n=0}^{+\infty} \frac{([-\boldsymbol{i}_\xi,-])^n}{n!}\ d-d= e^{[-\boldsymbol{i}_\xi,-]}d- d= e^{-\boldsymbol{i}_\xi} d e^{\boldsymbol{i}_\xi}-d. \]
The second equality follows from:
\[ e^{-\boldsymbol{i}_\xi}*0 = d\boldsymbol{i}_\xi -\frac{1}{2}[\boldsymbol{i}_\xi, d\boldsymbol{i}_\xi] + \frac{1}{3!}[\boldsymbol{i}_\xi,[\boldsymbol{i}_\xi, d\boldsymbol{i}_\xi]] + \ldots= \]
\[= \l_\xi - \boldsymbol{i}_{d\xi} - \frac{1}{2}[\boldsymbol{i}_\xi, \l_\xi]+ \frac{1}{2}[\boldsymbol{i}_\xi, \boldsymbol{i}_{d\xi}] +\frac{1}{3!}[\boldsymbol{i}_\xi,[\boldsymbol{i}_\xi, \l_\xi]] -\frac{1}{3!}[\boldsymbol{i}_\xi,[\boldsymbol{i}_\xi, \boldsymbol{i}_{d\xi}]]+ \ldots  = \l_\xi - \boldsymbol{i}_{d\xi} - \frac{1}{2}\boldsymbol{i}_{[\xi,\xi]}= \l_\xi,\]
where we used Lemma \ref{lemma eq cartan homotopy}, the definition of the holomorphic Lie derivative and  the fact that $\xi \in \MC_{\KS_X}(A)$.  \end{proof}
Now we can describe the sheaves of $\xi$-holomorphic differential forms: 
\[ \Omega^1_\xi =\{ \omega \in \A_\xi^{1,0} \mid d\omega \in \A_\xi^{2,0} \} = \{ \omega \in \A_\xi^{1,0} \mid e^{\boldsymbol{i}_\xi} (d+\l_\xi) e^{-\boldsymbol{i}_\xi} \omega \in \A^{2,0}_\xi\} =    \]
\[= \{ \omega \in \A_\xi^{1,0} \mid  (d+\l_\xi) e^{-\boldsymbol{i}_\xi} \omega \in \A^{2,0}_X\otimes A \}=\{ \omega \in \A_\xi^{1,0} \mid (\deltabar+\l_\xi) e^{-\boldsymbol{i}_\xi} \omega =0 \};  \]
in general:
\[ \Omega^p_\xi =\{ \omega \in \A_\xi^{p,0} \mid d\omega \in \A_\xi^{p,0} \}= \{ \omega \in \A_\xi^{p,0} \mid (\deltabar+\l_\xi) e^{-\boldsymbol{i}_\xi} \omega =0   \}.\]
\begin{remark}
The analysis of deformations of complex structures done in this section is coherent with the one done in Example \ref{ex.def varieta}. 
Infact, given a Maurer-Cartan element $\xi \in A^{0,*}_X(\T_X) \otimes \mathfrak{m}_A$, its associated deformation of complex structure of $X$ is $\A^{0,1}_\xi$, that defines a complex manifold whose sheaf of holomorphic functions is
\[ \Oh_\xi=\{f\in \A_X^0 \mid (\deltabar +\l_\xi) f=0\}= \ker (\deltabar + \l_x: \A_X^{0,0}\otimes A \to \A_X^{0,1}\otimes A);\] exactly as done by $\Phi$ in Example \ref{ex.def varieta}.   
\end{remark}

\section{Infinitesimal deformations of Higgs bundles} \label{sec.def Higgs}
This section is dedicated to the introduction of deformations of Higgs bundles and of pairs (manifold, Higgs bundle) and to the description of the results we have obtained on these problems; we will relate complete proofs in Section \ref{Sect.proof of main theorem}.

\begin{definition}
Let $X$ be a compact complex manifold. A \emph{Higgs bundle} on $X$ is a pair $(E,\theta)$, where $E$ is a holomorphic vector bundle on $X$ and $\theta \in H^0(X, \Eps nd(E)\otimes \Omega_X^1)$, such that $\theta \wedge \theta=0$.  
\end{definition}

\begin{definition} \label{def. deformationi di (X,E,theta)}
Let $X$ be a compact complex manifold and let $(E,\theta)$ be a Higgs bundle on $X$. Let $A$ be a local Artinian $\C$-algebra. An \emph{infinitesimal deformation} of $(X,E,\theta)$ over $A$ is the data $(X_A, E_A,\theta_A)$, where:
\begin{itemize}
\item[-] $X_A$ is a deformation of $X$ over $A$, i.e. $X_A$ is a scheme with a cartesian diagram 
$$\xymatrix{ X \ar[d]\ar[r]^-{p} & X_A \ar[d]^-{\pi} \\
           \Spec\C \ar[r] & \Spec A, }$$ 
where $\pi$ is flat,
\item [-] $E_A$ is a locally free sheaf of $\Oh_{X_A}$-modules on $X_A$, flat over $A$, with a morphism $q: E_A \to E$, such that $q: E_A \otimes_A \C \to E$ is an isomorphism, 
	\item [-] $\theta_A \in H^0(X_A, \Eps nd (E_A)\otimes \Omega^1_{X_A|A})$, with $\theta_A \wedge \theta_A=0$ and the maps $p$ and $q$ transform  $\theta_A$ in $\theta$ on the closed point.
\end{itemize}
If the deformation $X_A$ is trivial, i.e. $X_A \cong X \times \Spec A$, the pair $(E_A, \theta_A)$ is an \emph{infinitesimal deformation} of the Higgs bundle $(E,\theta)$ over $A$.
\end{definition}

\begin{definition} \label{def.iso tra defs}
Let $X$ be a compact complex manifold and let $(E,\theta)$ be a Higgs bundle on $X$. Let $(X_A,E_A,\theta_A)$ and $(X'_A,E'_A,\theta'_A)$ be two deformations of $(X,E,\theta)$ over the local Artinian $\C$-algebra $A$. They are \emph{isomorphic}, if there exists a couple $(\phi, \psi)$, where:
\begin{itemize}
\item[-] $\phi: \mathcal O_{X_A} \to \mathcal O_{X'_A}$ is an isomorphism of sheaves, such that $\phi\otimes_A \C$ is the identity on $\Oh_X$, 
\item[-] $\psi: E_A\to E'_A$ is an isomorphism of sheaves of $\mathcal O_{X_A}$-modules, where the structure of $\mathcal O_{X_A}$-sheaf on $E'_A$ is given by $\phi$, such that $\psi \otimes_A \C$ is the identity on $E$,
\item[-] $(\phi,\psi)$ map $\theta_A$ on $\theta'_A$.   
\end{itemize}
Two deformations of the Higgs bundle $(E,\theta)$ are said to be \emph{isomorphic}, if there exists a couple of isomorphisms $(\phi,\psi)$ as above, with $\phi=\Id_{X\times \Spec A}$.
\end{definition}

\noindent The above definitions lead to introduce two \emph{deformation functors}: 
\[ \Def_{(X,E,\theta)}: \bf{Art}_{\C} \to \bf{Set},\] 
that associates to every local Artinian $\C$-algebra $A$ the set $\Def_{(X,E,\theta)}(A)$ of isomorphism classes of deformations of $(X,E,\theta)$ over $A$ and 
\[ \Def_{(E,\theta)}: \bf{Art}_{\C} \to \bf{Set}, \]
that associates to every local Artinian $\C$-algebra $A$ the set $\Def_{(E,\theta)}(A)$ of isomorphism classes of deformations of the Higgs bundle $(E,\theta)$ over $A$. 

Our main goal is to study these deformation functors finding out dglas that govern them. 
Fix a Higgs bundle $(E, \theta)$ on a complex manifold $X$, consider the graded vector space 
\[  \bigoplus_{p+q=*} A_X^{q,p}(\D^1(E)), \]
we define a dgla structure on it, taking as bracket  
\[ [\omega\otimes P, \eta\otimes Q]= \omega\wedge \eta \otimes[P,Q] + \omega\wedge \l_{\sigma(P)}(\eta) \otimes Q - (-1)^{\deg \omega\cdot \deg \eta}  \l_{\sigma(Q)}(\omega)\wedge \eta \otimes P, \]
and as differential
\[ d(\omega \otimes P)=\deltabar\omega \otimes P + [\theta, \omega\otimes P],\]
for all $\omega, \eta \in A_X^{*,*}$ and $P,Q\in \D^1(E)$, this structure is a modification of the one defined in Example \ref{Ex.End(E),D(E) dgla}.  
Consider now the  sub graded vector spaces of it:
\[ N^*= \bigoplus_{p+q=*}A^{0,p}_X(\Eps nd(E) \otimes \Omega_X^q) \quad \mbox{and} \]
\[ L^*= \bigoplus_{p+q=*,\ p<* } A^{0,p}_X(\Eps nd (E) \otimes \Omega_X^q) \oplus A^{0,*}_X(\D^1(E)),\]
define on them structures of dglas as sub dglas of $A_X^{*,*}(\D^1(E))$, noting that they are closed under differential and bracket. By definition, they enter in the following exact sequence of dglas:
\[ 0 \to N \to L \stackrel{\sigma}{\rightarrow} A_X^{0,*}(\T_X) \to 0. \]

\smallskip

Our main result is to prove that the dgla $L$ governs infinitesimal deformations of $(X,E,\theta)$, as a consequence, we obtain that the dgla $N$ governs infinitesimal deformations of $(E,\theta)$.
To prove these statements, in Propositions \ref{prop.def mappa} and \ref{prop.mappa sui def}, we will construct a natural transformation of deformation functors
\[ \Phi_L: \Def_{L} \to \Def_{(X,E,\theta)}   \]
and, in Propositions \ref{prop.iniettivita} and \ref{prop.suriettivita}, we will prove it is an isomorphism. 
Note that, the dg vector space $L$ is the total complex of the Dolbeault resolution of the complex of sheaves \[{\mathcal K}: \qquad 0 \to \D^1(E) \stackrel{[-,\theta]}{\longrightarrow} \Eps nd (E) \otimes \Omega_X^1 \stackrel{[-,\theta]}{\longrightarrow} \Eps nd (E) \otimes \Omega_X^2 \stackrel{[-,\theta]}{\longrightarrow}  \Eps nd (E) \otimes \Omega_X^3 \to \ldots, \] 
where differentials are defined using the brackets of the dgla $L$. Then, as easy consequence of dglas tecniques, we obtain a description of tangent and obstruction spaces for the functor $\Def_{(X,E,\theta)}$. Summing up, we have the following
\begin{theorem} \label{Teo.dgla controlla higgs bundle}
The dgla $L$ governs infinitesimal deformations of $(X,E,\theta)$. In particular the space of first order deformations of $(X,E,\theta)$ is canonically isomorphic to the first hypercohomology space of the complex of sheaves ${\mathcal K}$ 
and obstructions are contained in the second hypercohomology space of it.
\end{theorem}
Our construction of the transformation $\Phi_L: \Def_{L} \to \Def_{(X,E,\theta)}$ will restrict to an isomorphism of deformation functors 
\[\Phi_N: \Def_{N} \to \Def_{(E,\theta)},\]
since the dg vector space $N$ is the total complex of the Dolbeault resolution of the complex of sheaves  
\[ {\mathcal K}': \qquad 0 \to \Eps nd (E) \stackrel{[-,\theta]}{\longrightarrow} \Eps nd (E) \otimes \Omega_X^1 \stackrel{[-,\theta]}{\longrightarrow} \Eps nd (E) \otimes \Omega_X^2\stackrel{[-,\theta]}{\longrightarrow}  \Eps nd (E) \otimes \Omega_X^3 \to \ldots,\]
dglas theory gives results on tangent and ostruction spaces for $\Def_{(E,\theta)}$. We have the following:
\begin{corollary} \label{Cor.dgla Higgs bundle}
The dgla $N$ governs infinitesimal deformations of the Higgs bundle $(E,\theta)$. In particular the space of first order deformations of $(E,\theta)$ is canonically isomorphic to the first hypercohomology space of the complex of sheaves $\mathcal K'$ and obstructions are contained in the second hypercohomology space of it.
\end{corollary}
Infinitesimal deformations of Higgs bundles were studied by Biswas. In \cite[Theorem 2.5]{Biswas}, he identified first order deformations of the pair (manifold, Higgs bundle) with the first hypercohomology space of the complex $\mathcal K$ and, in \cite[Remark 2.8]{Biswas}, he observed that the second hypercohomology space of the trucation at the third term of the complex $\mathcal K'$ contains obstructions to deformations of a Higgs bundle.
Similar results for infinitesimal deformations of principal $G$-bundles can  be found in \cite[Theorems 2.3 and 3.1]{Biswas-Ramanan}.

\section{Proof of Theorem \ref{Teo.dgla controlla higgs bundle}} \label{Sect.proof of main theorem}
This section is devoted to the proofs of our results stated in Section \ref{sec.def Higgs}. Let $(E, \theta)$ be a Higgs bundle on a complex  manifold $X$, let $L$ and $N$ be the dglas introduced in the previous section; we prove that the dgla $L$ governs infinitesimal deformations of $(X,E,\theta)$ and, as a particular case, we obtain that $N$ governs deformations of $(E,\theta)$. 

\begin{proposition} \label{prop.def mappa}
The map 
\[ \Phi_L: \MC_{L} \longrightarrow \Def_{(X,E,\theta)},\]
given, for $A \in \bf{Art}_{\C}$ and for $(x,y) \in \MC_L(A)$, by the isomorphism class of the deformation $(X_A(\sigma(x)), E_A(x), \theta_A(y))$, where:
\begin{itemize}
\item[-] ${\mathcal O}_{X_A(\sigma(x))}=\mathcal O_A(\sigma(x))=\ker(\deltabar+\l_{\sigma(x)}: \A_X^{0,0}\otimes A\longrightarrow \A_X^{0,1}\otimes A)$, where $\l$ is the holomorphic Lie derivative,
\item[-] $E_A(x)= \ker(\deltabar+x: \A_X^{0,0}(E)\otimes A \longrightarrow \A_X^{0,1}(E)\otimes A)$, 
\item[-] $\theta_A(y)=e^{\boldsymbol{i}_{\sigma(x)}}(\theta+y) \in A_{X_A(\sigma(x))|A}^{1,0}(\Eps nd (E))$, where $\boldsymbol{i}$ is the contraction,
\end{itemize}\end{proposition}
\begin{proof}
We have to verify that, $(X_A(\sigma(x)),E_A(x),\theta_A(y))$ defined above is a deformation of $(X,E,\theta)$ over $A$.

The element $x$ satisfies the Maurer-Cartan equation in the dgla $A^{0,*}_X(\D^1(E))$, then, by Remark \ref{rm.mc vs gauge banale}, it is locally gauge equivalent to zero, i.e. there exist an open covering $\U=\{U_\beta\}_\alpha$ of $X$ and elements $a_\alpha \in A^{0,0}_X(D^1(\Eps))|_{U_\alpha}\otimes \mathfrak m_A$, such that $e^{a_\alpha}*x|_{U_\alpha}=0$ and, taking the principal symbol, $e^{\sigma(a_\alpha)}*\sigma(x)|_{U_\alpha}=0$. 
As we will explicitly compute in proof of Proposition \ref{prop.mappa sui def}, it follows from definition of gauge action that 
$e^{a_\alpha}: E_A(x)|_{U_\alpha} = \ker(\deltabar+x)|_{U_\alpha} \rightarrow \ker\deltabar|_{U_\alpha}=E|_{U_\alpha} \otimes A$ and $e^{\sigma(a_\alpha)}: \Oh_A(\sigma(x))|_{U_\alpha} = \ker(\deltabar+\l_{\sigma(x)})|_{U_\alpha}\rightarrow \ker\deltabar|_{U_\alpha}=\Oh_X|_{U_\alpha}\otimes A$ are isomorphisms. By the last isomorphism $\Oh_A(\sigma(x))$ is a $A$-flat sheaf and  $\Oh_A(\sigma(x))\otimes_A \C \cong \Oh_X$, while from the first one we deduce that $E_A(x)$ is $A$-flat and that $E_A(x)\otimes_A \C\cong E$.
Moreover, the composition $e^{\sigma(a_\alpha)} \circ e^{a_\alpha}: E_A(x)|_{U_\alpha} \to E|_{U_\alpha}\otimes A \cong {\Oh_X}^{\rm{rk} E} \otimes A \to {\Oh_A(\sigma(x))}^{\rm{rk} E}$ assures that $E_A(x)$ is a locally free sheaf $\Oh_A(\sigma)$ modules.

It remains to prove that $\theta_A(y)$ is a section in $H^0(X_A(\sigma(x)),\Eps nd (E_A(x)) \otimes \Omega^1_{X_A(\sigma(x))|A})$, such that $\theta_A(y) \wedge \theta_A(y)=0$.

Consider the $A_X^{0,1}(\Eps nd (E)\otimes \Omega_X^1)$-component of the Maurer-Cartan equation for the element $(x,y)$:
\[ 0= \left.\left[ d(x,y) +\frac{1}{2}[(x,y),(x,y)] \right]\right|_{A_X^{0,1}(\Eps nd (E)\otimes \Omega_X^1)}= \deltabar y + [\theta, x] + \frac{1}{2}([x,y]-(-1)^{\deg x\cdot \deg y} [x,y])=\]
\begin{equation}\label{eq.sez}= [\deltabar, y]-(-1)^{\deg{x}\cdot \deg{\theta}}[x,\theta]+ [x,y] + [\deltabar, \theta]= [\deltabar, y+\theta]+ [x,y+\theta]=[\deltabar+x, y+\theta],  \end{equation}
that is the action of the differential $\deltabar + x$ on the element $\theta+ y$, then 
\[ (\deltabar + x)(\theta+ y)=0 \qquad \Rightarrow \qquad((\deltabar + x)e^{-\boldsymbol{i}_{\sigma(x)}})\ \theta_A(y)=0.\]
We claim that the last equation is equivalent to the condition that $\theta_A(y)$ is a section in 
\[ H^0(X_A(\sigma(x)),\Eps nd (E_A(x)) \otimes \Omega^1_{X_A(\sigma(x))|A}).\]
Infact, recall that $E_A= \ker (\deltabar + x: \A_X^{0,0}(E) \otimes A \to \A_X^{0,1}(E) \otimes A)$ and that $\Omega_{X_A(\sigma(x))|A}^1= \ker ((\deltabar + \l_{\sigma(x)}) e^{-\boldsymbol{i}_{\sigma(x)}}: \A_{X_A(\sigma(x))|A}^{1,0} \to \A_X^{1,1} \otimes A)$, therefore 
\begin{equation} \label{eq due fasci} \xymatrix{   \Eps nd (E_A) \otimes \Omega^1_{X_A(\sigma(x))|A} \subseteq \ker(\A_{X_A(\sigma(x))|A}^{1,0}(\Eps nd (E)) \ar[rr]^-{(\deltabar + x) e^{-\boldsymbol{i}_{\sigma(x)}}}& & \A_{X}^{1,1}(\Eps nd (E))\otimes A).} \end{equation}
Since the deformations of $E$, $\Omega^1_{X}$ and $E\otimes \Omega_X^1$ involved are locally trivial, they preserve the ranks, thus the two sheaves in (\ref{eq due fasci}) have the same rank and coincide.     

To obtain the equation $\theta_A(y) \wedge \theta_A(y)=0$, consider the $A_X^{0,0}(\Eps nd (E)\otimes \Omega_X^2)$-component of the Maurer-Cartan equation for the element $(x,y)$:
\begin{equation}\label{eq.wedge} 0= \left.\left[ d(x,y) +\frac{1}{2}[(x,y),(x,y)] \right]\right|_{A_X^{0,0}(\Eps nd (E)\otimes \Omega_X^2)}= [\theta, y]+ \frac{1}{2}[y,y]=\end{equation}
\[ = \theta \wedge y -(-1)^{\bar{y}} y\wedge \theta + \frac{1}{2} (y\wedge y -(-1)^{\deg{y}\cdot \deg{y}}y\wedge y)= \theta \wedge y + y\wedge \theta + y\wedge y=  \]
\[ =(\theta+ y) \wedge (\theta +y),\]
where we used that $\theta \wedge \theta =0$. 
By Lemma \ref{Lemma.e iso}, $e^{\boldsymbol{i}_{\sigma(x)}}: \A^*_X \to \A^*_X$ 
is an isomorphism of graded algebras, then the above equality assures that $\theta_A(y) \wedge \theta_A(y)= 0$. 
\end{proof}  

\begin{proposition} \label{prop.mappa sui def}
The map $\Phi_L$ induces a natural transformation of deformation functors
\[ \Phi_L: {\Def}_{L} \longrightarrow \Def_{(X,E,\theta)}.\]
\end{proposition}
\begin{proof}
Let $(x,y), (x',y') \in \MC_L(A)$ be gauge equivalent and
let $a\in A_X^{0,0}(\D^1(E))\otimes \mathfrak{m}_A$ such that 
\[ e^a*(x,y)=(x',y') \qquad  \Rightarrow \qquad   e^a*x=x' \qquad  \Rightarrow\qquad  e^{\sigma(a)}*\sigma(x)=\sigma(x').\]
From the gauge action, we will obtain that
$\deltabar+ e^{\sigma(a)}*l_{\sigma(x)}=e^{\sigma(a)}\circ(\deltabar+\l_{\sigma(x)})\circ e^{-\sigma(a)}$, and in particular $e^{\sigma(a)}$ gives an isomorphism between $\Oh_A(\sigma(x))=\ker (\deltabar+\l_{\sigma(x)})$ and $\Oh_A(\sigma(x'))=\ker (\deltabar+ e^{\sigma(a)}*l_{\sigma(x)})$. Infact:
\begin{eqnarray*}  
\l_{\sigma(x')}  & = & e^{\sigma(a)}*\l_{\sigma(x)}  =   \l_{\sigma(x)}+\sum_{n=0}^{+\infty} \frac{([\sigma(a),-])^n}{(n+1)!}([\sigma(a),\l_{\sigma(x)}]-d\sigma(a))= \\
   & = & \l_{\sigma(x)}+\sum_{n=0}^{+\infty} \frac{([\sigma(a),-])^n}{(n+1)!}([\sigma(a),\l_{\sigma(x)}]+[\sigma(a),\deltabar])= \l_{\sigma(x)}+\sum_{n=1}^{+\infty} \frac{([\sigma(a),-])^n}{n!}(\deltabar+\l_{\sigma(x)})= \\ 
   & = &  \sum_{n=0}^{+\infty} \frac{([\sigma(a),-])^n}{n!}(\deltabar+\l_{\sigma(x)})-\deltabar= e^{[\sigma(a),-]}(\deltabar+\l_{\sigma(x)})-\deltabar = e^{\sigma(a)} \circ(\deltabar+\l_{\sigma(x)})\circ e^{-\sigma(a)}-\deltabar . 
\end{eqnarray*}  
Making calculation on the gauge relation $e^a*x=x'$, we obtain that $\deltabar+ e^{a}*x=e^{a}\circ(\deltabar+x)\circ e^{-a}$ and in particular $e^{a}$ is an isomorphism between $E_A(x)= \ker(\deltabar+x)$ and  $E_A(x')=\ker(\deltabar+x')$. Infact: 
\begin{eqnarray}  
x'  & = & e^{a}*x  =   x+\sum_{n=0}^{+\infty} \frac{([a,-])^n}{(n+1)!}([a,x]-da)= \nonumber \\
   & = & x+\sum_{n=0}^{+\infty} \frac{([a,-])^n}{(n+1)!}([a,x]+[a,\deltabar])= x+\sum_{n=1}^{+\infty} \frac{([a,-])^n}{n!}(\deltabar+x)=\nonumber\\ 
   & = &  \sum_{n=0}^{+\infty} \frac{([a,-])^n}{n!}(\deltabar+x)-\deltabar= e^{[a,-]}(\deltabar+x)-\deltabar = e^a \circ(\deltabar+x)\circ e^{-a}-\deltabar \nonumber . 
\end{eqnarray} 
Moreover $e^{a}$ is an isomorphism of sheaves of $\Oh_{X_A}$-modules, where the structure of sheaf of $\Oh_{X_A}$-modules on $\Eps'_A$ is the one induced by the isomorphism $e^{\sigma(a_i)}$.
 
It remains to prove that the isomorphisms $(e^{\sigma(a)}, e^a)$ map $\theta_A(y)$ in $\theta_A(y')$.
Consider the $A_X^{0,0}(\Eps nd(E) \otimes \Omega_X^1)$-component of the gauge equation: 
\begin{equation}\label{eq gauge} y'= \left.(e^a*y)\right|_{A_X^{0,0}(\Eps nd (E) \otimes \Omega_X^1)}= \left. \left[ y + \sum_{n=0}^{+\infty} \frac{[a,-]^n}{(n+1)!} ([a,y]- da)\right] \right|_{A_X^{0,0}(\Eps nd (E) \otimes \Omega_X^1)}= \end{equation}
\[=y + \sum_{n=0}^{+\infty} \frac{[a,-]^n}{(n+1)!} ([a,y]- [\theta, a])= y + \sum_{n=0}^{+\infty} \frac{[a,-]^n}{(n+1)!} ([a,y]+ [a,\theta])=        \]
\[ =y + \sum_{n=0}^{+\infty} \frac{[a,-]^n}{(n+1)!} ([a,y+\theta])= y + \sum_{n=1}^{+\infty} \frac{[a,-]^n}{n!} (y+\theta)= y + \sum_{n=0}^{+\infty} \frac{[a,-]^n}{n!} (y+\theta)- \theta,\]
then
\[ y'+\theta = e^{[a,-]}(y+\theta) \quad \Longrightarrow \quad y'+\theta = e^a \circ(y+\theta)\circ e^{-\sigma(a)}.\] 
Thus, composing with the isomorphisms given by the contractions, we obtain:
\[ e^{\boldsymbol{i}_{\sigma(x')}}(y'+\theta) = e^{\boldsymbol{i}_{\sigma(x')}} e^a e^{-\boldsymbol{i}_{\sigma(x)}} e^{\boldsymbol{i}_{\sigma(x)}}(y+\theta) e^{-\sigma(a)} \Rightarrow \theta_A(y') = e^{\boldsymbol{i}_{\sigma(x')}} e^a  e^{-\boldsymbol{i}_{\sigma(x)}} \theta_A(y) e^{-\sigma(a)}\] 
and the following diagram commutes:
\begin{equation} \label{eq corrispondenza sezioni}\xymatrix{ X_A(\sigma(x')) \ar[r]^-{\theta_A(y')} \ar[d]^-{e^{-\sigma(a)}}  & A^{1,0}_{X_A(\sigma(x'))|A}(\Eps nd (E))  \ar[r]^-{e^{-\boldsymbol{i}_{\sigma(x')}}} & A^{1,0}_X(\Eps nd (E))\otimes A \\
 X_A(\sigma(x)) \ar[r]^-{\theta_A(y)}  & A^{1,0}_{X_A(\sigma(x))|A}(\Eps nd (E))  \ar[r]^-{e^{-\boldsymbol{i}_{\sigma(x)}}} & A^{1,0}_X(\Eps nd (E))\otimes A \ar[u]^-{e^a}.             }
\end{equation}
This explicitly gives the required correspondence between $\theta_A(y)$ and $ \theta_A(y')$ via the isomorphisms $(e^{\sigma(a)}, e^a)$. 
\end{proof}

Now we want to prove that the natural transformation $\Phi_L$ is an isomorphism of deformation functors.
For the injectivity we need the following Lemmas, in which we use the above notations. 
\begin{lemma} \label{lemma.caso varieta}
Let $x, y \in \A_X^{0,1}(\mathcal T_X)\otimes \mathfrak m_A$ Maurer-Cartan elements. Let $F,G: \bf{Art}_{\C} \to \bf{Set}$ be the following functors
\[ F(A)=\{\ \mbox{isomorphisms of complexes }\  e^s : (\A_X^{0,*}
\otimes A ,\deltabar +\mathfrak l_x) \to (\A_X^{0,*} \otimes A,\deltabar
+\mathfrak l_y), \]
\[ \mbox{with} \ s \in \A_X^{0,0}(\mathcal T_X)\otimes \mathfrak m_A,\ \mbox{that specialize
to identity}\ \} \]
\[
G(A)=\{\ \mbox{isomorphisms of sheaves of $A$-module }\ \psi:\mathcal O_A(x) \to \mathcal O_A(y),\ \mbox{that specialize  to  identity} \ \}, \]
for all $A\in \bf{Art}_{\C}$. 
Then the restriction  morphism  $F \to G$ is surjective.
\end{lemma}
\begin{proof} See \cite [Lemma II.7.2]{Dona.Tesi}.\end{proof}
\begin{lemma} \label{lemma.caso fibrati}
Let $x,y \in \A_X^{0,1}(\mathcal E nd (E))\otimes \mathfrak m_A$ Maurer-Cartan elements. 
Let $F,G: \bf{Art}_{\C} \to \bf{Set}$ be the following functors
\[ F(A)=\{\ \mbox{isomorphisms of complexes }\ e^s : (\A_X^{0,*}(E)
\otimes A ,\deltabar +x) \to (\A_X^{0,*}(E) \otimes A,\deltabar
+y), \]
\[ \mbox{with}\ s \in \A_X^{0,0}(\Eps nd (E))\otimes \mathfrak m_A, \ \mbox{that specialize
to identity}\ \} \]
\[
G(A)=\{\ \mbox{isomorphisms of sheaves of $\mathcal O_X$-module }\ \psi:E_A(x) \to E_A(y),\ \mbox{that specialize  to  identity} \ \}, \]
for all $A\in \bf{Art}_{\C}$. 
Then the restriction  morphism  $F \to G$ is surjective.
\end{lemma}
\begin{proof} The proof is essentially the same as in \cite [Lemma II.7.2]{Dona.Tesi}; it is sufficient to substitute the Kodaira-Spencer dgla with the dgla $A_X^{0,*}(\Eps nd (E))$.
\end{proof}

\begin{lemma} \label{lemma.caso coppia}
Let $x, y \in \A_X^{0,1}(\D^1(E))\otimes \mathfrak m_A$ be Maurer-Cartan elements. 
Let $F,G: \bf{Art}_\C \to \bf{Set}$ be the following functors
\[ F(A)=\{\ \mbox{couples of compatible isomorphisms of complexes }\ (e^s, e^{\sigma (s)}),\ \mbox{where} \] 
\[ e^{\sigma(s)}:(\A_X^{0,*}\otimes A ,\deltabar +\mathfrak l_{\sigma(x)}) \to (\A_X^{0,*} \otimes A,\deltabar
+\mathfrak l_{\sigma(y)}) \ \mbox{and} \  e^s : (\A_X^{0,*}(E)\otimes A ,\deltabar +x) \to (\A_X^{0,*}(E) \otimes A,\deltabar
+y),\] 
\[\mbox{with}\ s \in \A_X^{0,0}(\D^1(E))\otimes \mathfrak m_A, \ \mbox{that specialize
to identity}\ \} \]
\[ G(A)=\{\ \mbox{couples of compatible isomorphisms of complexes }\ (\phi,\psi),\ \mbox{where} \] 
\[ \phi:{\mathcal{O}}_A(\sigma(x)) \to {\mathcal{O}}_A(\sigma(y)) \ \mbox{and} \  \psi: E_A(x)\to E_A(y), \ \mbox{that specialize
to identity}\ \},\]
for all $A\in \bf{Art}_{\C}$. 
The natural morphism $F\to G$ is surjective.
\end{lemma}
\begin{proof}
Let $(\phi,\psi)$ be a couple of compatible morphisms in $G(A)$. By Lemma \ref{lemma.caso varieta}, there exists $a\in A^{0,0}_X({\mathcal{T}}_X)\otimes \mathfrak m_A$, such that the isomorphism $e^a: (A^{0,*}_X\otimes A, \deltabar+ {\mathfrak l}_{\sigma(x)}) \to (A^{0,*}_X\otimes A, \deltabar+ {\mathfrak l}_{\sigma(y)})$ lifts $\phi$. 
Any lifting $b \in A^{0,0}_X(\D^1(E))\otimes \mathfrak m_A$ of $a$ defines an isomorphism $e^b: E_A(x) \to e^b(E_A(x))$ compatible with $e^a:{\mathcal O}_A(\sigma(x)) \to {\mathcal{O}}_A(\sigma(y))$. 
The composition $\psi \circ e^{-b}: e^b(E_A(x)) \to E_A(x) \to E_A(y)$ is compatible with 
$\phi \circ e^{-a}=\Id_{{\mathcal{O}}_A(\sigma(y))}$ and, by Lemma \ref{lemma.caso fibrati}, there exists $c \in A^{0,0}_X(\Eps nd (E))\otimes \mathfrak m_A$ such that the isomorphism $e^c: (A^{0,*}_X(E) \otimes A, \deltabar+ (e^{b}*x)) \to (A^{0,*}_X(E) \otimes A, \deltabar+ y)$ lifts $\psi \circ e^{-b}$. Then the couple $(e^a, e^{c}\circ e^b=e^{c\bullet b})$ defines compatible isomorphisms in $F(A)$, that lift $(\phi,\psi)$. Here we indicate with $\bullet$ the Baker-Campbell-Hausdorff product (\cite{Man Roma}), then $c\bullet b \in A_X^{0,0}(\D^1(E))\otimes \mathfrak m_A$.         
\end{proof}

Now we are ready to prove the following

\begin{proposition} \label{prop.iniettivita}
The map $\Phi_L: {\Def}_{L} \longrightarrow \Def_{(X,E,\theta)}$ is injective.
\end{proposition}

\begin{proof}
Let $(x,y),(x',y')\in \MC_L(A)$, such that the associated deformations $(X_A(\sigma(x)),E_A(x),\theta_A(y))$ and  $(X_A(\sigma(x')),E_A(x'),\theta_A(y'))$ are isomorphic. Let $\phi:{\mathcal{O}}_A(\sigma(x)) \to {\mathcal{O}}_A(\sigma(x'))$ and $\psi: E_A(x) \to E_A(x')$ be compatible isomorphisms of sheaves, that map $\theta_A(y)$ on $\theta_A(y')$. 
By Lemma \ref{lemma.caso coppia}, there exists $a\in A_X^{0,0}(\D^1(E)) \otimes \mathfrak m_A$, such that the compatible isomorphisms of complexes
\[ (\A^{0,*}_X    \otimes A, \deltabar + \l_{\sigma(x)} ) \stackrel{ e^{\sigma(a)} }{\longrightarrow}(\A^{0,*}_X \otimes A, \deltabar + \l_{\sigma(y)}) \quad \mbox{and}\] 
\[ (\A^{0,*}_X(E) \otimes A, \deltabar + x) \stackrel{e^a}{\longrightarrow} (\A^{0,*}_X(E) \otimes A, \deltabar + y), \]
lift the couple $(\phi, \psi)$, i.e. such that $e^a*x=x'$, as done in proof of Proposition \ref{prop.mappa sui def}.
By hypothesis, the isomorphisms $(\phi=e^{\sigma(a)}, \psi=e^a)$ map 
$\theta_A(y)$ on $\theta_A(y')$ in the sense of diagram (\ref{eq corrispondenza sezioni}), then 
$\theta_A(y')=e^{\boldsymbol{i}_{\sigma(x')}} e^a e^{-\boldsymbol{i}_{\sigma(x)}}\theta_A(y)e^{-\sigma(a)}$, i.e. $y' + \theta= e^a \circ (y+ \theta) \circ e^{-\sigma(a)}$, that is equivalent to the gauge equation $e^a*y=y'$, as we done in (\ref{eq gauge}). 
Putting together the two gauge equations, we obtain that $e^a*(x,y)=(x',y')$.
\end{proof}

Let's now prove the sujectivity of the transformation $\Phi_L$. At first we need to recall the following facts. 
\begin{remark} \label{Rm.sur caso coppia}
As a consequence of Propositions \ref{prop.def mappa} and \ref{prop.mappa sui def}, the natural transformation  $\Psi:\Def_{A_X^{0,*}(D^1(E))}\to\Def_{(X,E)}$ defined in Example \ref{ex.def (X,E)} is well defined on deformation functors and, by Lemma \ref{lemma.caso coppia}, it is injective. 
It is classically known that first order deformations of the pair $(X,E)$ are identified with the space $H^1(X,\D^1(E))$ and that obstructions to deformations are contained in $H^2(X,\D^1(E))$. On the other hand, as stated in Remark \ref{rm.def dgla}, dglas tecniques assure that the same is true for the functor $\Def_{A_X^{0,*}(D^1(E))}$. It is quite easy to prove that the natural transformation $\Psi:\Def_{A_X^{0,*}(D^1(E))}\to\Def_{(X,E)}$ induces an isomorphism on tangent spaces and an injective morphism on obstructions, then, by \cite[Proposition 2.17]{Man Pisa}, $\Psi$ is smooth and surjective. See \cite{Tesi.mia} for all details. 
\end{remark}

\begin{proposition} \label{prop.suriettivita}
The map $\Phi_L: {\Def}_{L} \longrightarrow \Def_{(X,E,\theta)}$ is surjective.
\end{proposition}
\begin{proof}
Let $(X_A,E_A,\theta_A)$ be a deformation of $(X,E,\theta)$ over $A$. By surjectivity of $\Psi:\Def_{A_X^{0,*}(D^1(\Eps))}\to\Def_{(X,\Eps)}$, stated in Remark \ref{Rm.sur caso coppia}, there exists $x\in \MC_{A_X^{0,*}(\D^1(E))}(A)$, such that $(X_A,E_A)=(X_A(\sigma(x)), E_A(x))$. 
Consider the element $y=e^{-\boldsymbol{i}_{\sigma(x)}}(\theta_A)-\theta \in A_X^{1,0}(\Eps nd (E))\otimes A$, and prove that the couple $(x,y)$ satisfies the Maurer-Cartan equation in the dgla $L$. Infact the $A_X^{0,2}(\D^1(E))\otimes \mathfrak m_A$-component of the Maurer-Cartan equation is zero, because $x\in \MC_{A_X^{0,*}(\D^1(E))}(A)$. The component in $A_X^{0,1}(\Eps nd(E) \otimes \Omega^1_X)\otimes \mathfrak m_A$ is 
\[  \deltabar y+ [\theta, x]+[x,y]=(\deltabar+x)(\theta +y)=((\deltabar+x)e^{-\boldsymbol{i}_{\sigma(x)}})\ \theta_A=0\]
as in (\ref{eq.sez}) and because $\theta_A \in H^0(X_A,\Eps nd (E_A)\otimes \Omega_{X_A|A})$.
The last component is in $A_X^{0,0}(\Eps nd(E) \otimes \Omega^2_X)\otimes \mathfrak m_A$: 
\[ [\theta, y]+\frac{1}{2}[y,y]= (\theta+y)\wedge (\theta+y) = e^{-\boldsymbol{i}_{\sigma(x)}}(\theta_A\wedge \theta_A)=0,  \]by calculation in (\ref{eq.wedge}). Moreover it obvious that $\Phi_L(x,y)=(X_A,E_A,\theta_A)$. 
\end{proof}

Summing up, the natural transformation $\Phi_L: {\Def}_{L} \to \Def_{(X,E,\theta)}$ is an isomorphism of deformation functors and the dgla $L$ governs infinitesimal deformations of $(X,E,\theta)$.  

Consider now the natural transformation:
$$\begin{array}{rrllr}
\Phi_N: &\Def_N(A)& \longrightarrow &\Def_{(E,\theta)}(A), & \qquad \forall A \in \bf{A}rt_{\C},  \\
& (x,y)& \longrightarrow & (\ker(\deltabar+x), \theta+y),&
\end{array}$$
that is the same as the transformation $\Phi_L$ in the semplified case in which differential operators are substituted by endomorphisms. Our calculation for $\Phi_L$ assures that $\Phi_N$ is an isomorphism and the dgla $N$ governs infinitesimal deformations of $(E,\theta)$.

\section{Infinitesimal deformations of Hitchin pairs} \label{subsection.Hitchin pair}
In this section we define Hitchin pairs and study their deformations, following results of Section \ref{Sect.proof of main theorem}. Our aim is to identify a dgla that governs these deformations and obtain a description of tangent and obstruction spaces.

\begin{definition}
Let $X$ be a compact complex manifold and let $L$ be a vector bundle on $X$. A \emph{Hitchin pair} on $X$  consists of a holomorphic vector bundle $E$ on $X$ and a section $\theta \in  H^0(X,\Eps nd(E)\otimes L)$, such that $\theta\wedge \theta=0$ as a section in $H^0(X, \Eps nd(E) \otimes (L \wedge L))$. For $L = \Omega^1_X$, we recover the definition of Higgs bundles.
\end{definition}

There is an obvious definition of \emph{infinitesimal deformation} of a Hitchin pair $(E, L ,\theta)$ over a local Artinian $\C$-algebra $A$: it is the data of a deformation $E_A$ of the vector bundle $E$ over $A$, with a section  $\theta_A \in H^0(X\times \Spec A, \Eps nd(E_A)\otimes L\otimes A)$, that deforms $\theta$ over $A$, as in Definition \ref{def. deformationi di (X,E,theta)}. Observe that in deformations we are interested in $X$ and $L$ deform trivially.
As in Definition \ref{def.iso tra defs}, there is an obvious notion of isomorphism of deformations of a Hitchin pair.

This leads to introduce the \emph{functor of deformations} of the Hitchin pair $(E,L,\theta)$:
\[ \Def_{(E,L, \theta)}: \bf{Art}_\C \to \bf{Set},\]
which associates to every local Artinian $\C$-algebra $A$ the set of isomorphism classes of deformations of $(E,L, \theta)$ over $A$. 

Fix now a Hitchin pair $(E,L,\theta)$ on $X$; the tensor product $\Eps nd (E) \otimes \bigwedge^* L$ has a natural structure of Lie algebra on $\C$:
\[ [\phi\otimes h, \psi\otimes l]= \phi\circ \psi \otimes h \wedge l -(-1)^{\deg h \cdot \deg l} \psi \circ \phi \otimes l\wedge h =[\phi, \psi]\otimes h\wedge l,   \] 
if we define on it the differential:
\[ d(\psi \otimes l)=[\theta, \psi\otimes l],   \]
we obtain a dgla-structure on $\Eps nd (E) \otimes \bigwedge^*L$.
Consider now the tensor product of the dg algebra $A_X^{0,*}$ of differential forms on $X$ with the dgla $\Eps nd (E) \otimes \bigwedge^*L$; following Example \ref{Ex.dgla prod tensoriale}, it can be endowed with a dgla structure.
Denote by $M$ this dgla, explicitly:
\[  M^*= \bigoplus_{p+q=*}A_X^{0,p}(\Eps nd (E) \otimes \bigwedge^qL), \]
the differential and the brackets are:
\[ dx=d(\omega \otimes f)=\deltabar\omega \otimes f +(-1)^{\deg \omega} \omega \otimes [\theta, f] \] 
\[[x,x']=[\omega \otimes f, \omega' \otimes f']= (-1)^{\deg \omega' \cdot \deg f} \omega\wedge \omega' \otimes [f,f'],  \]
for all $\omega, \omega' \in A_X^{0,*}$ and $f, f'\in \Eps nd (E) \otimes\bigwedge^*L$.

\begin{remark}
It is easy to verify that the dgla structure on $A_X^{0,*}(\Eps nd(E) \otimes \Omega_X^*)$ defined in Section \ref{sec.def Higgs} is the same as the dgla structure one can obtain substituting in the above definitions the vector bundle $L$ with the shaef of holomorphic $1$-forms $\Omega_X^1$.
But observe that the above construction can not be used to define a dgla structure on $A_X^{0,*}(\D^1(E) \otimes \Omega_X^*)$ because the sheaf of Lie algebras $\D^1(E)$ is not $\Oh_X$-linear.
\end{remark}

Taking into account our constructions and calculations of Section \ref{Sect.proof of main theorem}, consider the natural transformation:
$$\begin{array}{rrllr}
\Phi_M: &\Def_M(A)& \longrightarrow &\Def_{(E,L,\theta)}(A), & \qquad \forall A \in \bf{A}rt_{\C},  \\
& (x,y)& \longrightarrow & (\ker(\deltabar+x), \theta+y) &
\end{array}$$
defined as transformation $\Phi_L$ of Proposition \ref{prop.def mappa}, for the case in which differential operators are simply endomorphisms and $L$ replaces $\Omega_X$. Our calculations on $\Phi_L$ assures that $\Phi_M$ is an isomorphism. Note that the dg vector space $M$ is the total complex of the Dolbeault resolution of the complex of sheaves \[ {\mathcal C:} \qquad 0 \to \Eps nd (E) \stackrel{[-,\theta]}{\longrightarrow} \Eps nd (E) \otimes L \stackrel{[-,\theta]}{\longrightarrow} \Eps nd (E) \otimes (L\wedge L) \stackrel{[-,\theta]}{\longrightarrow}  \Eps nd (E) \otimes (L\wedge L\wedge L) \to \ldots,\]  
where differentials are defined using brackets of the dgla $M$;  dglas theory gives a description of tangent and ostruction spaces of $\Def_{(E,L,\theta)}$. Summing up, we have proved the following

\begin{theorem} \label{Teo.dgla L Higgs bundle}
The dgla $M$ governs infinitesimal deformations of the Hitchin pair $(E,L,\theta)$. In particular the space of first order deformations of $(E,L,\theta)$ is canonically isomorphic to the first hypercohomology space of the complex of sheaves ${\mathcal C}$ 
and obstructions are contained in the second hypercohomology space of it.
\end{theorem}

\section{$L_\infty$-algebras in deformation theory}
In this section we recall some basic aspects of $L_\infty$-algebras theory and its link with deformation theory. We mainly follow \cite{Man Roma}.   
\begin{definition}
A \emph{graded coalgebra} is a graded vector space $C=\bigoplus_{i\in \Z} C_i$ with a morphism of graded vector spaces $\Delta: C\to C\otimes C$ called coproduct.\\
The graded coalgebra $(C,\Delta)$ is coassociative if $(\Delta\otimes \Id) \otimes \Delta= (\Id \otimes \Delta)\otimes \Delta$ and it is cocommutative if $\Delta= T\Delta$, where $T$ is given by $T(v\otimes w)= (-1)^{\deg v \deg w} w\otimes v$.  \\
Let $(C,\Delta_C)$ and $(D, \Delta_D)$ be two graded coalgebras, a degree zero linear morphism $f: C\to D$ is a \emph{morphism of coalgebras} if $(f\otimes f)\Delta_C = \Delta_D f$.  
\end{definition}
\begin{definition}
Let $(C, \Delta)$ be a graded coalgebra. A \emph{coderivation} of degree $n$ on it is a linear map of degree $n$, $d\in \Hom^n(C,C)$, 
that satisfies the coLeibnitz rule
\[ \Delta d = (d \otimes \Id + \Id \otimes d) \Delta. \]
A coderivation d is called a \emph{codifferential} if $d \circ d = 0$.
\end{definition}

\begin{example}
Let $V$ be a $\Z$-graded vector space over $\K$. \\
The \emph{tensor coalgebra} generated by $V$ is defined to be the graded vector space
\[ T(V) = \bigoplus_{n=0}^{+\infty} \bigotimes ^n V   \]
endowed with the associative coproduct
\[ a(v_1\otimes \ldots \otimes v_n) = \sum_{k=1}^{n-1} (v_1\otimes \ldots \otimes v_k) \otimes (v_{k+1}\otimes \ldots \otimes v_n).\]
The \emph{reduced tensor coalgebra} generated by $V$ is the sub-coalgebra $\overline{T(V)}= \displaystyle\bigoplus_{n=1}^{+\infty} \bigotimes ^n V$.\\
Let $I$ be the homogeneous ideal of $T(V)$ generated by $\left\langle v\otimes w -(-1)^{\deg v \deg w}w\otimes v; \ \forall\ v, w\in V\right\rangle$.\\
The \emph{symmetric coalgebra} generated by $V$ is defined as the quotient 
\[S(V)=\bigoplus_{n=0}^{+\infty} \bigodot^n V, \quad \mbox{with} \quad  \bigodot^n V= \frac{\bigotimes^n V}{I\cap \bigotimes^n V},\]
endowed with the associative coproduct
\[\Delta(v_1\odot \ldots \odot v_n)= \sum_{k=1}^{n-1} \sum_{\sigma\in S(k,n-k) }\epsilon(\sigma) \ (v_{\sigma(1)}\odot \ldots \odot v_{\sigma(k)})\otimes (v_{\sigma(k+1)}\odot \ldots \odot v_{\sigma(n)}),\] 
where $S(k,n-k)$ indicates the set of permutations of $n$ elements, such that $\sigma(i)<\sigma(i+1)$, for all $i\neq k$, and $\epsilon(\sigma)=\pm 1$ is the sign determined by the relation in $\bigodot^n V$: $v_1\odot \ldots \odot v_n= \epsilon(\sigma)\  v_{\sigma(1)}\odot \ldots \odot v_{\sigma(n)}$. 
The \emph{reduced symmetric coalgebra} generated by $V$ is the sub-coalgebra of $S(V)$ given by $\overline{S(V)}= \displaystyle \bigoplus_{n=1}^{+\infty} \bigodot ^n V.$ Let $\pi: \overline{T(V)} \to \overline{S(V)}$ be the projection. 
\end{example}

For future use, a graded coalgebra $(C, \Delta)$ is called \emph{nilpotent}, if $\Delta^n = 0$ for $n >> 0$. 
It is \emph{locally nilpotent}, if it is the direct limit of nilpotent graded coalgebras or equivalently if $C=\cup_n \ker \Delta^n.$ The reduced tensor coalgebra and the reduced symmetric coalgebra generated by a $\Z$-graded vector space are locally nilpotent.

\begin{proposition}
Let $V$ be a graded vector space and let $(C,\Delta)$ be a locally nilpotent cocommutative graded coalgebra. The composition with the projection $p:\overline{S(V)}\to V$ defines a bijective map:
\[ \Hom(C, \overline{S(V)}) \stackrel{p\circ \ }{\longrightarrow} \Hom(C,V)   \] 
with inverse given by
\begin{equation}\label{formula morfismi}
 f \mapsto F= \sum_{n=1}^{+\infty} \frac{1}{n!}f^{\odot n}\circ \pi \circ \Delta^{n-1}.
\end{equation}
\end{proposition}
\begin{proof} See \cite{Man Roma}, Propositions VIII.18 and VIII.26. \end{proof}

\begin{proposition}
Let $V$ be a graded vector space and let $(C,\Delta)$ be a locally nilpotent cocommutative graded coalgebra. The composition with the projection $p:\overline{S(V)}\to V$ defines a bijective map:
\[ \Coder^n(C, \overline{S(V)}) \stackrel{p\circ \ }{\longrightarrow} \Hom^n(C,V)   \] 
with inverse given by
\begin{equation}\label{formula coderivazione1} 
q \mapsto Q = \pi \sum_{n=1}^{+\infty} \frac{1}{n!} (q\otimes \Id^{\otimes n}) \circ \Delta^n. 
\end{equation}
\end{proposition}
\begin{proof} See \cite{Man Roma}, Proposition VIII.33. \end{proof}
\begin{remark} 
Observe that, if $C= \overline{S(V)}$,
for $q= \sum_k q_k \in \Hom^n(\overline{S(V)}, V)$ and $v_1\odot \ldots \odot v_n \in \overline{S(V)}$, formula (\ref{formula coderivazione1}) gives:
\begin{equation} \label{formula coderivazione2}
Q(v_1\odot\ldots\odot v_n)=\sum_{k=1}^n\sum_{\sigma\in S(k,n-k)}
\varepsilon(\sigma)\ q_k(v_{\sigma(1)}\odot\ldots\odot v_{\sigma(k)})
\odot v_{\sigma(k+1)}\odot\ldots\odot v_{\sigma(n)}.
\end{equation}
\end{remark}
\begin{definition}
Let $V$ be a graded vector space, an $L_{\infty}$-\emph{structure} on $V$ is a sequence of linear maps of degree 1
\[ q_k\colon \bigodot^k V[1]\to V[1],\qquad \mbox{for}\ k\ge 1,\]
such that the coderivation $Q$ induced on the reduced symmetric coalgebra $\overline{S(V[1])}$ by the homomorphism $q =\sum_k q_k$ as in formula (\ref{formula coderivazione2}), is a codifferential.\\
An $L_{\infty}$-\emph{algebra} is indicated with $(V,q_i)$ and the morphisms $q_i$ are called the \emph{brackets} of the $L_\infty$-algebra.
\end{definition}
\begin{remark} \label{Rm.DGLA L-infinto}
A dgla $L$ has a natural structure of $L_{\infty}$-algebra. Infact, $L$ is a graded vector space and it can be verified that the brackets:
\[ q_1(x)=-dx, \qquad q_2(x\odot y)=(-1)^{\deg x}[x,y] \qquad \mbox{and} \qquad q_{k}=0, \mbox{\ \ for all \ } k\geq 3 \]
satisfy condition $Q\circ Q=0$. 
\end{remark}
\begin{definition}
Let $(V,q_i)$ and $(W,\hat{q_i})$ be two $L_{\infty}$-algebras, 
a \emph{morphism $f_\infty: (V,q_i) \to (W,\hat{q}_1)$ of $L_{\infty}$-algebras} is a sequence of linear maps of degree $0$ 
\[ f_k\colon \bigodot^k V[1]\to
W[1],\qquad \mbox{for}\ k\ge 1,\]
such that the morphism of coalgebras induced on the reduced symmetric coalgebras by $f=\sum_k f_k$, commutes with the codifferentials induced by the two $L_{\infty}$-structures of $V$ and $W$, i.e., with above notations, $F\circ Q= \hat{Q}\circ F$. \\
To verify that $f_\infty$ is an $L_\infty$-morphism it is sufficient to verify that
\begin{equation} \label{formula.morfismo alg Linfinito}
\sum_{a=1}^{+\infty} f_a \circ Q_n^a = \sum_{a=1}^{+\infty}\hat{q}_a\circ F_n^a, \qquad \mbox{for all}\ n \in \N.
\end{equation}
\end{definition}

\begin{remark}
Let $(V, q_i)$ be an $L_\infty$-algebra, the condition $Q\circ Q=0$ on the coderivation induced on the reduced simmetric coalgebra by the brackets $q_i$ implies that $q_1\circ q_1=0$, then $(V,q_1)$ is a differential complex.\\
Let $f_\infty:(V,q_i) \to (W,\hat{q}_i)$ be an $L_\infty$-morphism, its linear part $f_1$ 
satisfies the equation $f_1\circ q_1=\hat{q}_1\circ f_1$, thus $f_1$ is a morphism of differential complexes $(V,q_1)\to (W,\hat{q}_1)$ and it induces linear maps in cohomology $H^i(f_1): H^{i}(V) \to H^{i}(W)$.
\end{remark}

We are now ready to define a deformation functor associated to an $L_{\infty}$-algebra.

\begin{definition}
Let $(V,q_i)$ be an $L_{\infty}$-algebra, the \emph{deformation functor} associated to it is the functor $\Def_V:\mathbf{Art}_{\C}\to \mathbf{Set}$ defined, for all $A\in \bf{Art}_{\C}$, by:
\[  \Def_V(A)=\frac{\MC_V(A)}{\sim_{\rm homotopy}},\]
where \[\MC_{V}(A)=\left\{ x \in V[1]^0\otimes \mathfrak m_A \;\strut\left\vert\; \sum_{n\ge1}\frac{q_n( x ^{\odot n})}{n!}=0\right.\right\} \]
and the \emph{homotopy relation} is the following equivalent relation on $\MC_V(A)$: two elements $x, y \in \MC_V(A)$ are \emph{homotopy equivalent}, if there exists $z(t,dt)\in \MC_{V[t,dt]}(A)$, such that $z(0)=x$ and $z(1)=y$.
\end{definition}

\begin{remark} \label{Rm.Def DGLA L-infto}
Let $L$ be a dgla, as observed in Remark \ref{Rm.DGLA L-infinto}, it has an $L_{\infty}$-structure. The generalized Maurer-Cartan equation for the $L_{\infty}$-algebra $L$ is exactly the Maurer-Cartan equation for the dgla $L$, since $q_n=0$ for $n\geq 3$. Moreover, for a dgla the homotopy equivalence coincide with the gauge equivalence \cite[Corollary 7.4]{Fio-Man cone}. Then the deformation functors associated to $L$ as dgla and as $L_{\infty}$-algebra coincide.  
\end{remark}

Both $\MC$ and $\Def$ are functors from the category of $L_\infty$-algebras to the category of set-valued functors of Artin rings. The functor $\MC$ acts on an morphism $f_\infty: (V,q_i)\to (W, \hat{q}_i)$ of $L_\infty$-algebras in the following way: 
\[ \MC(f_{\infty}): \MC_V \to  \MC_W\]
is the natural transformation of functors given, for $A\in\bf{Art}_\C$ and $x \in \MC_V(A)$, by
\[ \MC_V(f_\infty)(x)= \sum_{n=1}^{+\infty} \frac{1}{n!}\ f_n (x^{\odot n});    \]
it preserves the homotopy equivalence and then it induces
a natural transformation of deformation functors:
\[ \Def(f_\infty): \Def_V \to \Def_W.\]
\begin{remark} \label{Rm.mappe di teo di ostruzione Linfinito}
Let $(V, q_i)$ be an $L_\infty$-algebra and $\Def_V$ the deformation functor associated to it.  
It can be proved that the tangent space to $\Def_V$ is the first cohomology space of the complex $(V,q_1)$, $H^1(V)$, and that obstructions are naturally contained in $H^2(V)$. \\
Let $f_\infty: (V, q_i)\to (W,\hat{q_i})$ be an $L_\infty$-morphism, the linear maps $H^1(f_1): H^1(V)\to H^1(W)$ and $H^2(f_1): H^2(V)\to H^2(W)$ are morphisms of tangent spaces and of obstruction spaces, respectively, compatible with the morphism $\Def(f_\infty): \Def_V \to \Def_W$.
\end{remark}

\section{The Hitchin map}
This section is devoted to the study of the Hitchin map via $L_\infty$-algebras deformation theory. We explicitate an $L_\infty$-morphism that induces the Hitchin map and we obtain a condition on obstructions to deform Hitchin pairs.

\begin{definition}
Let $X$ be a compact complex manifold and let $L$ be a vector bundle on $X$. For any Hitchin pair $(E,L,\theta)$ with $\rm{rk}E=r$, let $H(E,L, \theta) = (\tr(\theta), \ldots, \tr(\theta^r))$, where the products are done in the associative algebras $\Eps nd(E)\otimes \bigodot^kL$. This define the \emph{Hitchin map}:
\[ H: {\mathcal M}_{X,L,r} \to \bigoplus_{k=1}^{r} H^0(X, \bigodot^k L), \]
from the moduli space $\mathcal M_{X,L,r}$ of Hitchin pairs on $X$ of rank $r$ with fixed $L$, to the space of global sections of the vector bundles $\bigodot^k L$. 
\end{definition}

\begin{remark}
Other definitions of the Hitchin map can be found in litterature. It can be defined fixing an arbitrary base of the space of polynomial function on $r\times r$ matrices invariant under conjugation. All the obtained maps are linked by an automorphism of the codomain and they can be considered equivalent for our studies.  
\end{remark}
Let now describe the Hitchin map at infinitesimal level. 
Fix the Hitchin pair $(E,L, \theta)$ on the compact complex manifold $X$, the Hitchin map on infinitesimal deformations of $(E,L, \theta)$ over $A\in \bf{Art}_\C$ is given by:
\[ H: \Def_{(E,L, \theta)}(A) \to \bigoplus_{k=1}^{r} H^0(X, \bigodot^k L)\otimes \mathfrak m_A,\]
it associates to every deformation $(E_A, \theta_A) \in \Def_{(E,L, \theta)}(A)$ the sections   
$(\tr(\theta_A), \ldots, \tr(\theta_A^r))$, which deform sections $(\tr(\theta), \ldots, \tr(\theta^r))$ along the trivial deformations of $L$ over $A$.
In Theorem \ref{Teo.dgla L Higgs bundle}, we identify a dgla that controls deformations of the Hitchin pair $(E,L,\theta)$:
\[ \bigoplus_{p+q=*}A_X^{0,p}(\Eps nd(E) \otimes \bigwedge^q L).\]
On the other hand it is easy to verify that the functor of deformations of the sections $(\tr(\theta), \ldots, \tr(\theta^r))$ along trivial deformations of $L$ is isomorphic to the deformation functor associated to the abelian dgla $\bigoplus_{k=1}^r A_X^{0,*}(\bigodot^k L)[-1]$ via the natural transformation:  
$$\begin{array}{rllr}
\Def_{\bigoplus_{k=1}^r A_X^{0,*}(\bigodot^k L)[-1]} (A) & \longrightarrow & \Def_{(\tr(\theta), \ldots, \tr(\theta^r))}(A), & \qquad \forall A \in \bf{A}rt_{\C},  \\
(l_1,\ldots, l_r )& \longrightarrow & (\tr (\theta)+l_1, \ldots, \tr(\theta^r)+l_r). &
\end{array}$$
Taking into account these dglas interpretations of the functors involved, our aim is to obtain the Hitchin map as the map induced at the level of  deformation functors from an $L_\infty$-morphism:
\begin{equation}\label{la mappa che stiamo cercando} h:  \bigoplus_{p+q=*} A_X^{0,p}(\Eps nd(E)\otimes \bigwedge^qL)[1] \to\bigoplus_{k=1}^r A_X^{0,*}(\bigodot^k L).\end{equation}

\noindent  For $n\geq 1$ and $k \in \{1, \ldots , r \}$, define the linear degree zero maps:
\[ g^k_n:(A_X^{0,*}(\Eps nd(E) \otimes L)[1])^{\odot n}  \to A^{0,*}_X(\bigodot^k L),\]
given, for all $\omega_i\in A_X^{0,*}$ and $f_i \in \Eps nd(E) \otimes L$, by
\[ g^k_n((\omega_1\otimes f_1)\odot \ldots \odot (\omega_n\otimes f_n)) = \omega_1\wedge \ldots \wedge \omega_n \otimes \overline{g^k_n}(f_1\odot \ldots \odot f_n),\]
where $ \overline{g^k_n}(f_1\odot \ldots \odot f_n)$ is the coefficient of $t_1\cdots t_n$ in 
\[ \tr(\theta + t_1 f_1 +\ldots +t_n f_n)^k.\]

\begin{remark}
For $n=1$, the map $g^k_1: A^{0,*}_X(\Eps nd(E) \otimes L)[1] \to A_X^{0,*}(\bigodot^k L)$ is
the identity on forms and it is given by $g^k_1(f)=k \tr(f\theta^{k-1})$ on $f\in \Eps nd (E) \otimes L$. 
It induces a morphism of complexes: 
\[ g^k_1: \bigoplus_{p+q=*} A_X^{0,p}(\Eps nd(E)\otimes \bigwedge^qL)[1] \to A_X^{0,*}(\bigodot^k L);\]
infact the definitions of the differentials on the two complexes and the fact that $g^k_1$ is the identity of forms reduce the proof to the following equality:
\[ g^k_1([\theta,f])=k \tr([\theta,f]\theta^{k-1})=\tr(\theta f\theta^{k-1}-f\theta^{k})= 0, \qquad \mbox{for all} \ f \in \Eps nd (E) \otimes L.\]
\end{remark}

\begin{proposition} \label{prop.g Linfty morf}
For all $k \in \{1, \ldots r\}$, the map 
\[ g^k_\infty: \bigoplus_{p+q=*} A_X^{0,p}(\Eps nd(E)\otimes \bigwedge^qL) \to A_X^{0,*}(\bigodot^k L)[-1], \]
defined by $\{g^k_n\}_n$ is an $L_\infty$-morphism.
\end{proposition}
For our proof we need the following
\begin{lemma}\label{lemma.polinomio man}
Let $A,B$ matrices and $k>0$, the coefficient of $t$ in $\tr\left((A+t[B,A])^k\right)$ is zero. \end{lemma}
\begin{proof} We prove by induction on $k$ that the $t$ coefficient in $(A+t[B,A])^k$ is $[B,A^k]$.
For $k=1$, there is nothing to prove. For $k >1$, using induction, we have 
\[ [B,A^k]=[B,A]A^{k-1}+A[B,A^{k-1}]=\sum_{i=0}^{k-1}A^i[B,A]A^{k-1-i},\]
where the last term is exactly the $t$ coefficient in $(A+t[B,A])^k$.
\end{proof}
\begin{proof}[Proof of Proposition \ref{prop.g Linfty morf}]
By formula (\ref{formula.morfismo alg Linfinito}), the map $g^k$ is an $L_\infty$-morphism if and only if
\[ \sum_{a=1}^{+\infty} g_a^k\circ Q_n^a = \sum_{a=1}^{+\infty} q_a \circ (G^k)_n^{a}, \qquad \mbox{for all}\ n \in \N.  \]
Since $g^k_n$ vanish on all terms in $A_X^{0,p}(\Eps nd(E)\otimes \bigwedge^qL)$, for $q \neq 1$, we only have to prove the formula on $A_X^{0,*}(\Eps nd(E)\otimes L)^{\odot n}$ and on $A_X^{0,*}(\Eps nd(E))\odot A_X^{0,*}(\Eps nd(E)\otimes L)^{\odot n-1}$. 

Let's start with $y_i = \omega_i \otimes f_i \in A_X^{0,*}(\Eps nd(E)\otimes L)$. Since the only non zero brackets are $q_1$ and $q_2$ in the first $L_\infty$-algebras and $q_1$ in the second one, we have:
\[  \sum_{\sigma \in S(1, n-1)} \epsilon(\sigma)\ g_n^k\left( q_{1}(y_{\sigma(1)}) \odot y_{\sigma(2)}\odot\ldots \odot y_{\sigma(n)}\right) =\]\[=\sum_{\sigma \in S(1, n-1)}\epsilon(\sigma)\ \deltabar \omega_{\sigma(1)} \wedge \ldots \wedge \omega_{\sigma(n)} \cdot \overline{g_n^k}(f_1\odot \ldots \odot f_n)=\] 
\[ = \deltabar(\omega_1 \wedge\ldots \wedge \omega_n)\cdot \overline{g_n^k}(f_1\odot \ldots \odot f_n) = q_1(g^k_n(y_1 \odot \ldots \odot y_n))   \]
and the formula is proved. 

Now let  $x=\omega\otimes f \in A_X^{0,p}(\Eps nd(E))$ and $y_i =\omega_i \otimes f_i \in A_X^{0,p_i}(\Eps nd(E)\otimes L)$, the formula to be proved becomes:
\[ g^k_n(q_1(x)\odot y_{1} \odot \ldots \odot y_{n-1})+  \sum_{\sigma\in S(1,n-2)} \epsilon(\sigma)\ g^k_{n-1}\left( q_2(x\odot y_{\sigma(1)})\odot \ldots \odot y_{\sigma(n-1)}\right) =0,\]
Now observe that, for all  $x$ and $y_i$ as above:  
\[ q_1(x)= -\deltabar \omega \otimes f - (-1)^p \omega \otimes [\theta, f] \qquad \mbox{and} \qquad q_2(x\otimes y_i)= (-1)^p (\omega\wedge \omega_i) \otimes [f,f_i]. \]
Then the above formula becomes:
\[ (-1)^p(\omega \wedge \omega_1 \wedge \ldots \wedge \omega_{n-1}) \cdot \left( \overline{g^k_n}([f, \theta] \odot f_{1} \odot \ldots \odot f_{n-1})+ \right.\] 
\begin{equation} \label{eq coefficienti}
+ \sum_{\sigma\in S(1,n-2)} \left.\overline{g^k_{n-1}} ([f, f_{\sigma(1)}]\odot \ldots \odot f_{\sigma(n-1)})\right) =0,  \end{equation}
where the Koszul sign $\epsilon(\sigma)$ disappears when we order differential forms.
Consider now the power:
\[ \left(\theta +t_1 f_1 + \ldots + t_{n-1} f_{n-1}+ t[f,\theta +t_1 f_1 + \ldots + t_{n-1} f_{n-1}]\right)^k; \]
the trace of the $t\cdot t_1 \cdots t_{n-1}$ coefficient is exactly the left term in formula (\ref{eq coefficienti}) and Lemma \ref{lemma.polinomio man}, applied with $A=\theta +t_1 f_1 + \ldots + t_{n-1} f_{n-1}$ and $B=f$, assures it is zero. 
\end{proof}
As immediate consequence of the above proposition, we have the $L_\infty$-morphism (\ref{la mappa che stiamo cercando}) we are looking for:
\begin{proposition} \label{Prop.h L-infty morphism}
The map 
\[ h:  \bigoplus_{p+q=*} A_X^{0,p}(\Eps nd(E)\otimes \bigwedge^qL)[1] \to\bigoplus_{k=1}^r A_X^{0,*}(\bigodot^k L),\]
defined by $h=(g^1, \ldots , g^r)$ is an $L_\infty$-morphism.
\end{proposition}

To identify the natural transformation induced on deformation functors by the $L_\infty$-morphism $h$ with the Hitchin map, we need the following: 

\begin{lemma} \label{prop.formula di pola} 
With the above notations, for all $k\in \{ 1, \ldots ,r \}$, we have: 
\[ \sum_{n=0}^k\frac{1}{n!}\ g^k_n(y^{\odot n})= \tr((\theta+y)^k-\theta^k), \qquad \mbox{for all} \ y \in A_X^{0,0}(\Eps nd(E) \otimes L).\]
\end{lemma}
\begin{proof} It is the polarization formula and follows from an easy calculation.
\end{proof}

\begin{proposition} \label{Prop.h Hitchin map}
The natural transformation of deformation functors
\[ \Def(h)(A): \Def_{(E,L,\theta)}(A) \to \bigoplus_{k=1}^{r} H^0(X, \bigodot^k L)\otimes \mathfrak m_A, \qquad \forall \ A \in \bf{Art}_\C,\] 
induced by the $L_\infty$-morphism $h$ 
is the Hitchin map.
\end{proposition}
\begin{proof}
Let $(x,y)\in (A_X^{0,0}(\Eps nd (E)) \oplus A^{0,1}_X(\Eps nd (E) \otimes L)) \otimes \mathfrak m_A$ be a Maurer-Cartan element, with $A \in \bf{Art}_\C$; by Lemma \ref{prop.formula di pola}: 
\[ \Def(h)(x,y)= \left(\tr((\theta + y)^k - \theta^k)\right)_{k=1\ldots r}. \]
On the other hand, the image via the Hitchin map of the infinitesimal deformation $(E_A=\ker (\deltabar +x), \theta_A=\theta+y)$ associated to $(x,y)$ is 
\[ (\tr(\theta_A), \ldots , \tr(\theta_A^r))=\left(\tr (\theta_A^k)\right)_{k=1 \ldots r}. \]
To conclude observe that 
\[ \tr (\theta_A^k)= \tr (\theta+y)^k = \tr\theta^k + (\Def(h)(x,y))_k,\]
then the image of $(E_A,\theta_A)$ via the Hitchin map is the deformation of $(\tr(\theta^k))_{k=1\ldots r}$ given by the element $\Def(h)(x,y)$.  
\end{proof}
The above deformation theoretic interpretation of the Hitchin map leads to the following description of obstructions to deform Hitchin pairs.
\begin{corollary} \label{Cor.ostruzione e mappa Hitchin}
The obstructions to deform the Hitchin pair $(E,L,\theta)$ are contained in the kernel of the linear map
\[H^2(h_1): \mathbb{H}^2({\mathcal C}) \to \bigoplus_{k=1}^{r} H^1(X, \bigodot^k L), \] 
\[H^2(h_1)([\omega \otimes f])=[(k \cdot \omega\otimes \tr(f\theta^{k-1}))_{k=1\ldots r}], \qquad \mbox{for} \ \omega \otimes f \in A_X^{0,1}(\Eps nd (E) \otimes L).\]
\end{corollary}
\begin{proof}
We use a general argument on theory of deformation via $L_\infty$-algberas (Remark \ref{Rm.mappe di teo di ostruzione Linfinito}). Since the Hitchin map $H: \Def_{(E,L,\theta)} \to \Def_{(\tr(\theta), \ldots, \tr(\theta^r))}$ is induced by the $L_\infty$-morphism $h:\bigoplus_{p+q=*} A_X^{0,p}(\Eps nd(E)\otimes \bigwedge^qL)[1] \to\bigoplus_{k=1}^r A_X^{0,*}(\bigodot^k L)$ 
the linear map
\[ H^2(h_1) : \mathbb{H}^2({\mathcal C}) \to \bigoplus_{k=1}^{r} H^1(X, \bigodot^k L)\]
is a morphism of obstruction theories, i.e., it commutes with the natural obstruction
maps for $\Def_{(E,L,\theta)}$ and $\Def_{(\tr(\theta), \ldots, \tr(\theta^r))}$. Since the dgla $\bigoplus_{k=1}^r A^{0,*}_X(\bigodot^kL)$ is abelian, the last deformation functor is unobstructed and the obstructions to deform $(E,L,\theta)$ are annihilated by $H^2(h_1)$. 
\end{proof}


\begin{thebibliography}{99}
\begin{raggedright}
\bibitem{Biswas} I. Biswas, \emph{A remark on a deformation theory of Green and Lazarsfeld,} J. reine angew. Math. {\bf{449}}, (1994), 103-124.
\bibitem{Biswas-Gothen-Logares} I. Biswas, P. Gothen, M. Logares, \emph{On moduli space of Hitchin pairs}, \texttt{arXiv:math/09124615} 
\bibitem{Biswas-Ramanan} I. Biswas, S. Ramanan \emph{An infinitesimal study of the moduli of Hitchin pairs,} J. London Math. Soc. (2), no.2, (1994), 219-231.
\bibitem{Fio-Man cone}D. Fiorenza, M. Manetti: \emph{$L_{\infty}$-structures on mapping
cones,} Algebra \& Number Theory, \textbf{1}, (2007), 301-330.
\bibitem{Fio-Man periodi} D. Fiorenza, M. Manetti, \emph{$L_\infty$-algebras, Cartan homotopies and period maps,} \texttt{arXiv:math/0605297}.
\bibitem{Fio-Man periodi generalizzati} D. Fiorenza, M. Manetti, \emph{A period map for generalized deformations,} J. Noncommut. Geom. 3, (2009), 579-597.  
\bibitem{Goldman-Millson 1} W. Goldman, J. Millson, \emph{The deformation theory of representations of foundamental groups of compact K\"{a}hler manifolds,} Publ. Math. I.H.E.S., 67, (1988), 43-96.
\bibitem{Goldman-Millson 2} W. Goldman, J. Millson, \emph{The homotopy invariance of the Kuranishi space,} Illinois Journal of Math., 34, (1990), 337-367. 
\bibitem{Hitchin} N . J . Hitchin, \emph{The self-duality equations on a Riemann surface,} Proc. London Math. Soc. (3), 55, (1987), 59-126.
\bibitem{Kodaira} K. Kodaira, \emph{Complex manifold and deformation of complex structures,} Springer-Verlag (1986).
\bibitem{Dona.Tesi} D. Iacono, \emph{Differential graded Lie algebras and deformations of holomorphic maps,} PhD Thesis, \texttt{arXiv:math/0701091}. 
\bibitem{Man Pisa} M. Manetti, \emph{Deformation theory via differential graded Lie algebras,} Seminari di Geometria Algebrica 1998-1999, Scuola Normale Superiore (1999).
\bibitem{Man Roma} M. Manetti, \emph{Lectures on deformations of complex manifolds,} Rend. Mat. Appl., (7) \textbf{24}, (2004), 1-183. 
\bibitem{Man Kahler} M. Manetti, \emph{Cohomological constraint to deformations of compact K\"ahler manifolds,} Adv. Math., {\bf 186}, (2004), 125-142.
\bibitem{Tesi.mia} E. Martinengo, \emph{Higher bracket and moduli space of vector bundles,} PhD Thesis.
\bibitem{Sernesi} E. Sernesi, \emph{Deformation of algebraic schemes,} Springer, \textbf{334}, (2006). 
\bibitem{Simpson} C. T. Simpson, \emph{Higgs bundles and local systems,} Publ. Math. I.H.E.S. {\bf{75}}, (1992), 5-95.
\bibitem{Simpson 2} C. T. Simpson, \emph{Moduli of representations of the fundamental group of a smooth projective 
variety, II,} Publ. Math. I.H.E.S. {\bf{80}}, (1994), 5–79.
\bibitem{Voisin} C. Voisin, \emph{Th\'eorie de Hodge et g\'eom\'etrie alg\'ebrique complexe, I,} Soci\'et\'e Mathématique de France, Paris, (2002). 
\end{raggedright}
\end{thebibliography}
\end{document}